\documentclass[10pt,reqno]{amsart}
\usepackage{graphicx}
\usepackage{indentfirst,csquotes}

\usepackage{amssymb,amsthm,amsmath}
\usepackage{xcolor,paralist,hyperref,fancyhdr,etoolbox}

\theoremstyle{plain}

\newtheorem{theorem}{Theorem}[section]

\newtheorem{corollary}[theorem]{Corollary}

\newtheorem{lemma}[theorem]{Lemma}
\newtheorem{proposition}[theorem]{Proposition}

\theoremstyle{remark}

\newtheorem{definition}[theorem]{Definition}
\newtheorem{example}[theorem]{Example}

\newtheorem{remark}[theorem]{Remark}

\hypersetup{ colorlinks=true, breaklinks=true, linkcolor=blue, citecolor = [RGB]{193, 37, 37} }

\newcommand{\Prob}[2][]{\mathbb{P}{#1} \left[ #2 \right]}

\newcommand{\Esp}[2][]{\mathbb{E}{#1} \left[ #2 \right]}

\newcommand{\tr}[1]{\operatorname{tr} \left( #1 \right)}

\newcommand{\diag}[1]{\operatorname{diag} \left( #1 \right)}

\newcommand{\Abs}[1]{\left| #1 \right|}

\newcommand{\Norm}[1]{\left\lVert #1 \right\rVert}

\newcommand{\Braces}[1]{\left\lbrace #1 \right\rbrace}

\newcommand{\Brackets}[1]{\left[ #1 \right]}

\newcommand{\Parentheses}[1]{\left( #1 \right)}

\newcommand{\Angles}[1]{\left\langle #1 \right\rangle}

\newcommand{\Ind}[1]{\mathbf{1}_{ #1 } }

\newcommand{\diff}[0]{\operatorname{d}}

\newcommand{\Diff}[0]{\operatorname{D}}

\newcommand{\cadlag}{\emph{càdlàg}}

\newcommand{\EE}[0]{\mathbb{E}}
\newcommand{\FF}[0]{\mathbb{F}}

\newcommand{\HH}[0]{\mathbb{H}}

\newcommand{\NN}[0]{\mathbb{N}}
\newcommand{\PP}[0]{\mathbb{P}}

\newcommand{\RR}[0]{\mathbb{R}}

\newcommand{\Bb}[0]{\mathcal{B}}

\newcommand{\Dd}[0]{\mathcal{D}}
\newcommand{\Ee}[0]{\mathcal{E}}
\newcommand{\Ff}[0]{\mathcal{F}}

\newcommand{\Hh}[0]{\mathcal{H}}

\newcommand{\Ll}[0]{\mathcal{L}}
\newcommand{\Mm}[0]{\mathcal{M}}
\newcommand{\Nn}[0]{\mathcal{N}}

\newcommand{\Pp}[0]{\mathcal{P}}

\newcommand{\Rr}[0]{\mathcal{R}}
\newcommand{\Ss}[0]{\mathcal{S}}

\newcommand{\Xx}[0]{\mathcal{X}}

\begin{document}


\title[Coupled FBSDEs with jumps in random environments]{Coupled forward-backward stochastic differential equations with jumps in random environments} 
\date{\today}
\thanks{This work was supported by the National Council of Science and Technology (CONACyT), scholarship number 863210, as part of the PhD Thesis of the corresponding author.}


\author{Daniel~Hernández-Hernández}
\address{Centro de Investigación en Matemáticas A.C. Calle Jalisco s/n. 36240 Guanajuato, México}
\email{dher@cimat.mx}

\author{Joshué Helí~Ricalde-Guerrero}
\address{Centro de Investigación en Matemáticas A.C. Calle Jalisco s/n. 36240 Guanajuato, México}
\email{joshue.ricalde@cimat.mx}


\begin{abstract}
In this paper we obtain results for the existence and uniqueness of solutions to coupled Forward-Backward Stochastic Differential Equations (FBSDEs) with jumps defined on a random environment. This environment corresponds to a measured-valued process, similar to the one found in Conditional McKean-Vlasov Differential Equations and Mean-Field Games with Common Noise. The jump term in the FBSDE is dependent on the environment through a stochastic intensity process. We provide examples which relate our model with FBSDEs driven by Cox and Hawkes processes, as well as regime-switching Conditional McKean-Vlasov differential equations. 



\end{abstract} 


\maketitle



\section{Introduction}
\label{Sec:Introdution}
Let $(\Omega,\Ff,\PP)$ be a probability space equipped with a general filtration $\FF=\{\Ff_t\}$ satisfying the usual assumptions of right continuity and completeness, and $\Ff_T=\Ff$, where $T$ is a fixed finite horizon. Given an exogenous measured-valued process $\mu$, defined on this probability space and adapted to $\FF$, in this paper we are interested in solving a stochastic forward-backward system of the form 
 \begin{align}
    \label{Eq:Forward-Component}
    X_t^\mu
    =&
    X_0^\mu
    +
    \int_0^t b^\mu\Parentheses{ s, X_s^\mu, Y_s^\mu, Z_s^\mu, U_s^\mu(\ \cdot\ ) } \diff t
    +
    \int_0^t \sigma^\mu \Parentheses{ s, X_s^\mu, Y_s^\mu, Z_s^\mu, U_s^\mu(\ \cdot\ ) } \diff W_t
    \\ \nonumber
    &+
    \int_{ [0,t] \times \mathbf{R} } \gamma^\mu \Parentheses{
        s, r, X_{s-}^\mu, Y_{s-}^\mu, Z_{s-}^\mu, U_{s-}^\mu(\ \cdot\ )
    } \tilde{\eta}^\lambda( \diff s, \diff r ),
    \\ \label{Eq:Backward-Component}
    Y_t^\mu
    =&
    Y_T^\mu
    +
    \int_t^T f^\mu\Parentheses{ s, X_s^\mu, Y_s^\mu, Z_s^\mu, U_s^\mu(\ \cdot\ ) } \diff s
    -
    \int_t^T Z^\mu_s \diff W_s
    \\ \nonumber
    &-
    \int_{(t,T] \times \mathbf{R}} U_{s-}^\mu(r) \tilde{\eta}^\lambda( \diff s, \diff r )
    -
    (M^\mu_T - M^\mu_t),
    \\ \label{Eq:Initial-Terminal-Conditions}
    X_0^\mu \sim& \mu_0 \in \Pp_2(\RR^d; \Ff_0),
    \qquad
    Y_T^\mu = g(\mu_T,X_T) \in L^2(\RR^n; \Ff_T), \;\;M_0=0.
\end{align}
Here $W$ is a Wiener process, $\eta^\lambda$ is the lifting of a multivariate marked Poisson process on $[0,T] \times \mathbf{R}$ with predictable compensator $K^\lambda(t, \diff r) \diff t$, $\tilde{\eta}^\lambda$ is the corresponding compensated martingale, and $M$ is an $\RR^n$-valued, $\FF$-martingale orthogonal to $(W,\tilde{\eta}^\lambda)$.

The study of Forward-Backward Stochastic Differential Equations (FBSDEs) can be traced  back to the second half of the past century; see, for instance, the fundamental works \cite{bismut_conjugate_1973}, \cite{bismut_introductory_1978}, \cite{pardoux_adapted_1990}, \cite{peng_general_1990}, \cite{el_karoui_backward_1997}, \cite{peng_fully_1999}, \cite{ma_forward-backward_2007}.   These type of systems arise naturally in the study of problems in stochastic optimal control, in particular in the formulation of necessary optimality conditions for controlled processes via the  stochastic version of the Pontryagin's Maximum Principle (PMP) \cite{yong_stochastic_1999}. This class of BSDEs are explicitly dependent on the dynamics of the control, giving rise to  a coupled system of FBSDEs. 

The simpler  class of FBSDEs consists of that where the only source of randomness in the system is due to the noise introduced by a  Brownian motion. Sufficient conditions for the existence and uniqueness of solutions for fully-coupled FBSDEs were studied in \cite{peng_fully_1999}, where the authors made use of specific monotonicity conditions (usually referred  as $G$-monotonicity in the context of Functional Analysis \cite{oettli_existence_1998} \cite{brokate_generalized_nodate}), to guarantee the existence and uniqueness of solutions adapted to the filtration generated by the noise process. Their arguments are based on the \emph{method of continuation} \cite{yong_finding_1997}, which is  useful to solve FBSDEs derived from optimal control problems, and can also be modified to include noise processes of Poissoninan type \cite{tang_necessary_1994}, \cite{zhen_forward-backward_1999}.

Jump processes offer  one of the best suited tools to model atypical and unexpected events, and it is natural to include them in systems evolving according with a  FBSDEs. In fact, many interesting problems from queuing theory, operations research and population genetics require the manipulation of discontinuous dynamical systems, specifically those represented by integrals of Poisson random measures, see, for instance,  \cite{baccelli_elements_2002}, \cite{daley_introduction_2003}, \cite{daley_introduction_2008}, \cite{bremaud_changes_1978}. Doubly stochastic Poisson point processes and self-exciting Hawkes processes are often used in the field of finance and risk theory to model, for instance,  the behaviour of economic assets and options, the default risk of the banking markets and the dividends of assurance companies \cite{fouque_handbook_2013}, \cite{bacry_hawkes_2015}, \cite{hawkes_hawkes_2022}, \cite{hawkes_hawkes_2018}. For the study of optimal control problems for  jump-diffusion processes, we refer to  \cite{oksendal_applied_2019}, where the authors present an interesting introduction to the field.

In  the present work we include a   point processes $N^\lambda$, and the main novelty is that we assume that the corresponding intensity $\lambda$  depends on the exogenous measure-valued process $\mu$. This process acts as a random environment for the system due to the external source of randomness affecting explicitly  the forward-backward system. The intuition behind this phenomenon is the same found in the theory of  \emph{Mean-Field Games with Common Noise}. Namely, when a representative agent is interested in finding   an equilibrium when a large number of players are involved, an optimization problem naturally arises, subject to the dynamics of an external random flow of probability measures. In this context, FBSDEs influenced by a random environment appear in the formulation of the maximum principle, with $\mu$  taking values in the Polish space of Borel probability measures \cite{carmona_probabilistic_2018}.

The election of  jump processes depending on the environment has two main motivations. The first one comes from the need to obtain a sufficiently robust and  workable model that considers the effect of external inputs through periodic shocks. In many real-world applications, systems  are constantly affected by external factors, which do not necessarily evolve in continuous-time, and changes are observed after  certain positive random times, e.g. new arrivals to  queue systems or claims to assurance companies. In these cases, it is reasonable to consider  that periodic shocks are directly affected by the existing environment (arrivals diminish if there is a congestion before the entrance to the queue, and claims increase if the locality is exposed to natural disasters). By allowing the structure of the jumps process to have a relationship with the environment it is possible to incorporate exogenous information. This relationship is established in terms of the intensity of the jumps process, covering a large number of problems, while maintaining  meaningful results.
 
The second motivation comes from mean-field games theory. Succinctly, this theory is based on the fact that individual particles (players) participating in a highly populated dynamical system (game) are interchangeable with each other, and their individual actions become weaker as the population increases \cite{cardaliaguet_notes_nodate}, \cite{cardaliaguet_master_2019}, \cite{caines_large_2006}, \cite{carmona_probabilistic_2013}. Usually, the interactions only occur through some statistical information, like the empirical mean. 

The FBSDE derived from the maximum principle contains a measure-valued environment $\mu$ and,  when the game is at  equilibrium the environment coincides with the \emph{conditional law of the process itself}, resulting in an SDE of the \emph{conditional McKean-Vlasov type}, a kind of differential equation that has been studied since the 1970s in the theory of \emph{Propagation of Chaos} \cite{burkholder_topics_1991}, \cite{chaintron_propagation_2022}, \cite{chaintron_propagation_2022-1}. The relevant aspect here is that, owing to the information structure of the system (respectively, game), the particularities of $\mu$ (the conditional law of the player's state) are such that it is a process adapted to the filtration of the common noise. Historically, the study of these objects has been centered around SDEs in which the environment is adapted to a Gaussian-type of filtration, as it was originally presented in \cite{carmona_mean_2016}. To the best of our knowledge, the most general result in this direction  considers a Gaussian white noise field $W^0 = \Parentheses{ W^0( \Lambda, B ), \Lambda \in \Bb(\Xi), B \in \Bb(\RR_+) }$, see  section of \emph{Notes and Complements} in \cite[Chapter 2]{carmona_probabilistic_2018-1}). However, by directly embedding the information of the environment into the jump part, we add it to the type of possible common noise processes in MFGs and the FBSDEs that represent them, opening up the possibilities for a kind of McKean-Vlasov differential equation where the law is conditioned on a discontinuous, marked, Poisson-type of process.

A usual assumption in the study of FBSDEs is related with the information structure with which  the problem is posed, i.e. the filtration attached to the probability space is generated only by the noise processes influencing the evolution of the system; 
namely, the natural filtration of a Brownian motion and/or a Poisson point process. This critical assumption is due to the backward component of the equation and, within our framework,
the technical difficulties of including an external random environment   come from the fact that one has to work with a general filtration where
the martingale decomposition may no longer hold, due to  the presence of both the environment and the noise within the system. 
In this paper it is required that $\FF$  satisfies the usual conditions of completeness and, for the filtration generated by $(\mu,N^\lambda,W)$, it is assumed that it is  \emph{immersed} in $\FF$ \cite{bremaud_changes_1978}, \cite{burkholder_topics_1991}, \cite{azema_vershiks_2001}, which consists of guaranteeing the martingale property of a process  after some enlargement of the underlying filtration. In our context, thanks to the \emph{immersion condition} (or equivalently, the $\Hh$-hypotheses or the \emph{compatibility condition}, \cite{bremaud_changes_1978}, \cite{burkholder_topics_1991}, \cite{kurtz_weak_2014}, \cite{carmona_probabilistic_2018-1})  the martingale decomposition holds, which is reflected in the additional (orthogonal) martingale term $M^\mu$ in \eqref{Eq:Backward-Component}.
The approach presented in this paper to model the  jumps process, based on marked point processes with intensity functions depending on an exogenous random environment,  include a large class  of stochastic intensity processes. We describe  them as \emph{admissible (Poissonian) noise processes} and, naturally, they induce the notion of an \emph{admissible set-up} as any sufficiently rich probability space in which an admissible noise can be defined; see Theorem \ref{Thm:Existence-of-admissible-settings}.

To summarize, in  Section \ref{Prelim} we present the set of spaces for the stochastic processes involved in the FBSDE model, including the random environment. This model is presented  rigorously in Section \ref{MR}, where we establish  the existence and uniqueness of (strong) solutions to \eqref{Eq:Forward-Component}-\eqref{Eq:Initial-Terminal-Conditions} in Theorem \ref{Thm:Existence-Uniqueness-FBSDE}. For readability, we present preliminary results for  uncoupled systems in Section \ref{PUFBE}, while the proof of the main result is presented in Section \ref{Sec:FBSDEs}. Finally, these results are applied to several examples in Section \ref{Sec:SE}, including the regime switching.

\section{Notation and Basic Notions}
\label{Prelim}
Throughout this paper, we use $(x\cdot y) = x^\top y$ to denote the usual inner product in $\RR^d$. If $G$ is a given full rank matrix on $\RR^{n \times d}$, we write $(x\cdot y)_G := (Gx)^\top y$ for  $x \in \RR^d$ and  $y \in \RR^n$. We also use $c_G < \infty$ to denote the positive proportionality constant $( x \cdot Gx)_G = c_G (x\cdot x)$. For matrices, we use the Frobenius inner product, i.e. $(A \cdot B)_\mathrm{Fr} = \tr{A B^\top}$ for $A,B$ real-valued matrices of the same dimensions, with the corresponding norm given by $\Norm{A}_\mathrm{Fr}^2 = \tr{ A A^\top }$.

Given a non-empty time interval $I \subset \RR_+$ and a Polish space $(\Xx,d_\Xx)$, we denote by 
\begin{align*}
    \mathrm{D}(I;\Xx)
    :=
    \Braces{X:I \to \Xx \ \middle\vert \ X\mbox{ has \cadlag\ trajectories }},
\end{align*}
endowed with the Borel $\sigma$-algebra $\Dd(I;\Xx)=\Bb\Parentheses{ \mathrm{D}(I;\Xx) }$, with respect to the Skorohod topology \cite{billingsley_convergence_2013}.   Denote by $\Pp(\Xx)$ the set of Borel probability measures on $\Xx$ and,  for    $p \in [1,\infty)$, define the set of probability measures with finite moments of order $p$ as
$$
  \Pp_p(\Xx):=
            \Bigg\lbrace 
                \nu \in \Pp(\Xx) 
                \ :\ 
                \int_\Xx d_\Xx(x_0,x)^p \nu( \diff x ) < +\infty
                \ 
                \mbox{for  any } x_0 \in \Xx
            \Bigg\rbrace.
$$
Given  $\nu^1,\; \nu^2\in \Pp(\Xx) $, define $\Pi(\nu^1,\nu^2) \subset \Pp(\Xx \times \Xx)$ as the set of probability measures over the product space $\Xx \times \Xx$ with marginals $\nu^1$ and $\nu^2$, respectively. The $p$-th Wasserstein distance between two probability measures $\nu^1, \nu^2 \in \Pp_p(\Xx)$ is then defined as
\begin{equation}
    \label{Eq:W_p^p}
    W_p (\nu^1, \nu^2)^p
    :=
    \min \Braces{ 
        \int_{ \Xx \times \Xx } d_\Xx( x_1, x_2 )^p \pi( \diff x_1, \diff x_2 )
        \ :\ 
        \pi \in \Pi( \nu^1, \nu^2 )
    }.
\end{equation}
Given $\Theta$ a stochastic process, we introduce some filtrations associated with $\Theta$. The natural filtration is denoted by $\FF^{ \mathrm{nat},\Theta }$, while the  $\PP$-null  complete right-continuous augmentation  is denoted by  $\FF^{ \Theta }$. In general, we say that a filtration $\FF$ satisfies the usual conditions if it is right-continuous and complete. Product spaces will always be endowed with the product $\sigma-$field.
        
Any mapping $\Psi:\Omega \times [0,T] \times \Pp_2(\RR^d) \times \RR^d \to \RR^p$ is denoted by $\Psi_t(\nu,x)(\omega)$ for each  $(\nu,x) \in \Pp_2(\RR^d) \times \RR^d$ and $(\omega, t) \in \Omega \times [0,T]$. Moreover, for each $(\nu,x),(\nu',x')\in \Pp_2(\RR^d) \times \RR^d$, their difference is denoted by $ \delta \Psi_t(\nu,x)(\omega)= \Psi\Parentheses{ \omega, t, \nu, x } - \Psi\Parentheses{ \omega, t, \nu', x' }$. This notation is also used to represent the difference of any two elements $x,x'$ of an Euclidean space; i.e. $ \delta x = x - x'$. 

        {\bf Spaces for stochastic processes.} 
        We write $\Hh^2(\RR^d)$ (or simply $\Hh^2$ when there is no ambiguity) for the Hilbert
space of $\RR^d$-valued, square-integrable martingales (with respect to $\FF$), 
i.e.
\begin{align*}
	\Hh^2(\RR^d)
	:=
	\Braces{ 
		M:\Omega \times [0,T] \to \RR^d
		\ \middle|\ 
		M\mbox{ is an }\FF\mbox{-martingale with }
		\sup_{ [0,T] }\Esp{ \Abs{ X_t }^2 }<\infty
	},
\end{align*}
endowed with  the norm $\Norm{\ \cdot\ }_{\Hh^2(\RR^d)}$,
induced by  the inner product $\Parentheses{ M\cdot M' }_{\Hh^2}=\Esp{ M_T M'_T }$, defined  by $\Norm{ M }_{ \Hh^2 }:=\Esp{ \Abs{ M_T }^2 }^{ \frac{1}{2} }$.  For any subset $B \subset \Hh^2$, $B^\perp \subset \Hh^2$ denotes the orthogonal set of $B$, i.e. the set of martingales in $\Hh^2$ which are orthogonal to all martingales in $B$. The space of well-posed integrands with respect to $M$  is denoted by $\HH^2(M)$, and is endowed  with the  predictable quadratic variation norm $\Norm{\ \cdot\ }_{\HH^2(M)}$. If $M'$ is another square integrable martingale, we write $\HH^2(M,M') = \HH^2(M) \cap \HH^2(M')$. We denote by $\Ll^2(M)$ the space of martingale-driven, Itô stochastic integrals of processes in $\HH^2(M)$ with respect to $M$, i.e.  $\Ll^2(M):=\Braces{ (H \bullet M)\ :\ H \in \HH^2(M) },$ equipped with the  martingale norm \cite{protter_stochastic_2005}.  We use the notation $[H,H']$ to denote the \textit{quadratic covariation} between semimartingales $H$ and $H'$, and $\Angles{H,H'}$ to denote its compensator when it exists (known also as \textit{conditional} or \textit{predictable quadratic covariation} \cite[Ch. III.5]{protter_stochastic_2005}). In particular,  we write $[H]=[H,H]$ and $\Angles{H} = \Angles{H,H}$. The space of adapted, \cadlag\ processes with finite supremum norm is denoted by $(\Ss^2, \Norm{\ \cdot\ }_\ast$). We  use the symbol $\FF' \hookrightarrow \FF$ whenever $\FF'$ and $\FF$ are filtrations on a $\sigma$-algebra $\Ff$ such that $\FF' \subset \FF$ and  $\FF'$-martingales are also $\FF$-martingales. In particular, an $\FF$-adapted \cadlag\ process $\Theta = \Parentheses{ \Theta_t, 0 \leq t \leq T }$ with values in a Polish space is said to be \emph{compatible with $\FF$ (under $\PP$)} if $\FF^\Theta\hookrightarrow\FF$.

{\bf Random environment.}  The main element of the system (\ref{Eq:Forward-Component})-(\ref{Eq:Backward-Component}) consists of an auxiliary, square-integrable process $\mu$ taking values on  the 
Wasserstein space $(\Pp_2({\RR^d}),W_2)$. This process plays the role of an exogenous random 
environment influencing the system. More precisely, let 
$\mu$ be a measurable map from $(\Omega,\Ff,\PP)$ into 
\begin{align*}
	\Parentheses{ 
		\mathrm{D} \Parentheses{ [0,T]; \Pp_2({\RR^d}) }, 
		\Dd \Parentheses{ [0,T]; \Pp_2(\RR^d) } 
	}.
\end{align*}
For each  $0 \leq t \leq T$, denote by $\mu_t$  the projection of $\mu$ at time  $t$,
i.e.  $\mu_t(\omega) := e_t \circ \mu(\omega)$. We say that $\mu$ is a 
\emph{Random Environment }if there exists 
an $\Ff_0^{\mu}$-measurable variable $X_0$ and a probability measure 
$\nu_0 \in \Pp_2(\RR^d)$,  the \emph{initial distribution of $\mu$}, such that
\begin{align}
	\label{Eq:Conditions-on-Random-Environment-I}
	\Ll(X_0) &= \nu_0,
	&
	\Prob{  \mu_0 = \nu_0} &= 1,
\end{align}
and, for each $\nu\in\Pp_2(\RR^d)$,
\begin{align}
	\label{Eq:Conditions-on-Random-Environment-II}
	\Esp{ \sup_{[0,T]} W_2(\mu_t,\nu)^2 } < \infty.
\end{align}

We are interested in \emph{Marked Poisson processes with random intensity and marks on a measurable space $(\mathbf{R},\Rr)$, such that the information from the environment influences their behaviour}. This is done assuming, first, that $\FF^{\mu} \subset \FF$. Then, we say that a non-negative random function $K^\lambda : \Omega \times [0,T] \times \Rr \to \RR_+$ is an \emph{admissible $\FF$-intensity kernel} if it is an $\FF$-predictable and locally integrable kernel from $(\Omega \times [0,T],\Ff \otimes \Bb([0,T])$ into $(\mathbf{R},\Rr)$; i.e.,
\begin{itemize}
    \item 
        $K^\lambda(\cdot,\cdot,C)$ is $\FF$-predictable for each $C \in \Rr$,
    \item 
        $K^\lambda(\omega,t,\cdot)$ is a finite measure on $(R,\Rr)$ for each $(\omega,t) \in \Omega \times [0,T]$.
\end{itemize}
The \textit{Admissible Noise Process} is defined as a marked Poisson process $(N^\lambda,\xi)$ on $[0,T]$ with marks on $\mathbf{R}$, with $\FF$-intensity kernel $K^\lambda$; i.e., $N^\lambda$ is a Poisson process of stochastic intensity $\lambda_t(\omega) = K^\lambda(\omega,t,\mathbf{R})$ and jumps occur at times $\tau = \Braces{\tau_i, i \geq 1 }$. The second component $\xi = \Braces{ \xi_i, i \geq 1 }$ is a sequence of $\mathbf{R}$-valued random variables with transition probabilities
\begin{align*}
    \Prob{ \xi_n \in C \middle\vert \sigma\Braces{ (\tau_i, \xi_i), i = 1, \ldots, n-1 } \vee \sigma\Braces{ \tau_n } }
    =
    \frac{ K^\lambda(\omega, \tau_n, C ) }{ K^\lambda(\omega, \tau_n, \mathbf{R} ) }.
\end{align*}
Note that $(N^\lambda,\xi)$ can also be characterised by its so-called \textit{lifting}; that is, the random measure $\eta^\lambda$ defined as
\begin{align}\label{lifting}
    &\eta^\lambda(B) = \sum_{i \geq 1} \delta_{(\tau_i,\xi_i)}(B),
    &
    & B \in \Bb([0,T]) \otimes \Rr,
\end{align}
with predictable compensator $K^\lambda(\omega,t, \diff r) \diff t$.

In the unmarked case, we simply say $N^\lambda$ is an $\FF$-admissible noise if $(N^\lambda,\xi)$ is a process with marks on $\mathbf{R}=\RR$ and intensity kernel $K^\lambda(\omega,t, \diff r) = \lambda_t(\omega) \delta_{ \{1\} }(\diff r)$ for some non-negative process $\lambda$, locally integrable and predictable with respect to $\FF$. The process $N^\lambda$ can be built following \cite[Thm. 5.4.1]{bremaud_point_2020}; namely,  using a standard Poisson process with parameter one, introduce a change of time. Then, the process $N^\lambda$ is such that the compensated process 
\begin{align*}
    \tilde{N}^\lambda = \Parentheses{ N^\lambda_t - \int_0^t \lambda_s \diff s, 0 \leq t \leq T }
\end{align*}
is an $\FF-$ martingale. 

\begin{remark}
In general, when we refer to $(N^\lambda,\xi)$ as an $l$-dimensional process, we mean the family
\begin{align*}
    (N^\lambda,\xi) = \Braces{ (N^{\lambda(j)},\xi^{(j)}), j = 1, \ldots, l },
\end{align*}
where each $(N^{\lambda(j)},\xi^{(j)})$ is an independent Poisson process with an independent $\FF$-intensity kernel $K^{\lambda(j)}$ from $(\Omega \times [0,T], \Ff \otimes \Bb([0,T]))$ into $(\mathbf{R},\Rr)$. Similarly, we write $K^{\lambda}$ (resp. $\eta^\lambda$) for the vector with $j$-th entry $K^{\lambda(j)}$ ($\eta^{\lambda(j)}$).
\end{remark}

{\bf Admissible Set-up.} We can now describe an \emph{Admissible Set-up} as an space endowed with both an \emph{idiosyncratic} and  \emph{admissible noises}.
\begin{definition}
\label{Def:Admissible-Set-Up}
A tuple $(\Omega,\Ff,\FF,\PP,X_0,\mu,\lambda,(N^\lambda,\xi),W)$ is said to be an \emph{admissible set-up} for the problem if
\begin{enumerate}
    \item
	$\mu$ is a random environment on $(\Omega, \Ff, \PP)$;
    \item
	$X_0 \in L^2( \Omega, \Ff_0, \PP; \RR^d )$;
		
    \item
        $K^\lambda$ is a non-negative, $\FF$-predictable, locally integrable kernel from $(\Omega \times [0,T],\Ff \otimes \Bb([0,T])$ into some measurable space $(\mathbf{R},\Rr)$;
		
    \item
        $(N^\lambda,\xi)$ is an $l$-dimensional Marked Poisson process, adapted to $\FF$, with $\FF$-stochastic intensity kernel $K^\lambda$ and lifting $\eta^\lambda$;

    \item
        $W$ is a $k$-dimensional $\FF$-Wiener process, sometimes called the {\it idiosyncratic noise}.
		
    \item
        $(X_0, \mu, \eta^\lambda)$ is independent of $W$ (under $\PP$);
		
    \item
	Any $\FF^{(X_0, \mu, \eta^\lambda, W)}$-martingale is also an $\FF$-martingale (under $\PP$).
\end{enumerate}
\end{definition}

\begin{remark}
Condition 7 above is commonly known in literature as \emph{$\Hh$-hypothesis}, meaning that $\FF^{(X_0, \mu, \eta^\lambda, W)}\hookrightarrow\FF$, or simply that ${(X_0, \mu, \eta^\lambda, W)}$ is \emph{compatible} with $\FF$, see \cite{aksamit_enlargement_2017}.
\end{remark}

In order to describe properly  the solution to the system  (\ref{Eq:Forward-Component})-(\ref{Eq:Backward-Component}) we need to define the following space.
\begin{align}
    \nonumber
    \mathfrak{H}^\lambda
    :=
    \Bigg\lbrace
        U : \Omega \times [0,T] \times \mathbf{R} \to \RR^{n \times l}
		\Bigg|\ &
		U \text{ is } \FF \otimes \Rr-\text{predictable, with } 
        \Norm{ U(\omega,t,\cdot) }_{\lambda_t ( \omega )} < \infty,
        \\ \label{Eq:Frac-H-lambda}
		&\diff \PP \otimes \diff t\mathrm{-a.e. }
    \Bigg\rbrace,
\end{align}
where $\Norm{\ \cdot\ }_{\lambda_t (\omega)}^2 = \Parentheses{ \cdot, \cdot }_{\lambda_t (\omega)}$ is a random norm defined for matrix-valued functions, induced by the inner product
\begin{align}
    \label{Eq:lambda-inner-product}
    &\Parentheses{ u \cdot v }_{\lambda_t (\omega)} 
    := 
    \int_\mathbf{R} \tr{ u(r) \diag{K^\lambda(\omega,t;\diff r)} v(r)^\top },
    &
    &\diff \PP \otimes \diff t\mbox{-a.e. }(\omega, t) \in \Omega \times [0,T].
\end{align}
This is due to the fact that, for each $U \in \mathfrak{H}^\lambda$ (writing $U_t(\cdot) = U(t,\cdot))$, when
\begin{align}
    \label{Eq:Finite-Predictable-Quadratic-variation-N^lambda}
    &\int_0^T \Norm{U_t(\cdot)}_{\lambda_t}^2 \diff t < \infty
    &
    &\PP-\mbox{a.s.},
\end{align}
the integral process defined below is well-defined for $\diff \PP \otimes \diff t$-almost all $(\omega,t) \in \Omega \times [0,T]$:
\begin{align}
    \label{Eq:Integral-wrt-eta^lambda}
    \int_{[0,t] \times \mathbf{R} } U(\omega,s,r) {\eta}^\lambda (\omega, \diff s, \diff r)
    :=
    \sum_{j=1}^l \sum_{i \geq 1} U^{(\cdot,j)} (\omega,\tau_i^{(j)}(\omega),\xi_i^{(j)}(\omega)) \Ind{\tau_i^{(j)}(\omega) \leq t },
\end{align}
where $\tau^{(j)}$ (resp. $\xi^{(j)}$) is the sequence of jumping times (marks) of the $j$-member of $(N^\lambda,\xi)$, and $U^{(\cdot,j)}$ denotes the $j$-th column of $U$. Moreover, if the random variable in \eqref{Eq:Finite-Predictable-Quadratic-variation-N^lambda} is also integrable with respect to $\PP$, then the integrals with respect to the compensated martingale measure
\begin{align*}
    \tilde{\eta}^\lambda(\omega, \diff s, \diff r)
    :=
    {\eta}^\lambda(\omega, \diff s, \diff r)
    -
    K^{\lambda}(\omega, s, \diff r) \diff s
\end{align*}
are also  square-integrable martingales, see \cite[Thm. 5.1.33-34]{bremaud_point_2020}. 

We refer \eqref{Eq:Integral-wrt-eta^lambda} as the \emph{integral of $U$ with respect to $(N^\lambda,\xi)$} --equivalently, the \emph{integral of $U$ with respect to $\eta^\lambda$}--. Following the notation used for well-posed integrals with respect to square-integrable martingales, we write
\begin{align*}
    \HH^2(\tilde{\eta}^\lambda)
    :=
    \Braces{
        U \in \mathfrak{H}^\lambda
        \  \middle\vert\ 
        \Esp{ \int_0^T \Norm{U_t(\cdot)}_{\lambda_t}^2 } < \infty
    }
\end{align*}
for the space of well-posed integrands with respect to $(N^\lambda,\xi)$, and 
\begin{align*}
    \Ll^2(\tilde{\eta}^\lambda)
    :=
    \Braces{
        \Parentheses{
            \int_{[0,t] \times \mathbf{R} } U(s,r) \tilde{\eta}^\lambda (\diff s, \diff r)
        }_{t \in [0,T]}
        \ \middle\vert\ 
        U \in \HH^2(\tilde{\eta}^\lambda)
    }
\end{align*}
for the space of martingales that can be expressed as compensated integrals with respect to $(N^\lambda,\xi)$.

\begin{remark}
Observe that for the unmarked case, equation \eqref{Eq:lambda-inner-product} reduces to
$\Parentheses{ u \cdot v }_{\lambda_t} = \tr{ u \diag{ \lambda_t } v^\top }$ and $
	\int_0^t U_s \diff \tilde{N}^\lambda_s
    :=
    \sum_{j=1}^l \Big\lbrace
        \sum_{\tau_j \leq t} U_{\tau_j}^{(\cdot,j)} 
        - 
        \int_0^t U_s^{(\cdot,j)} \lambda^{(j)}_s \diff s
    \Big\rbrace,
$ for $t \in [0,T]$.
\end{remark}

\section{Main Results}
\label{MR}
We now proceed to describe precisely our main results in a filtered probability space $(\Omega,\Ff,\FF,\PP),$ where an admissible set-up shall be built. We adopt the usual approach to define  an $l$-dimensional Poisson process, which is given as a vector $N = (N^{(1)},\ldots,N^{(l)})$, where each entry is an independent copy from the rest, $Q^{(i)}$ is a probability measure on $(\mathbf{R},\Rr)$, and each coordinate
\begin{align*}
    &N^{(i)}([0,\cdot] \times C \times D)
    -
    \int_{ [0,\cdot] \times C \times D } \diff t \otimes   Q^{(i)}(\diff r) \otimes \diff s
\end{align*}
is an $\FF^N$-martingale with respect to a given reference measure $\PP^0$ (specified below), for all $C \in \Rr$ and $D \in \Bb(\RR_+)$ such that 
\begin{align*}
    &\int_{ [0,t_0] \times C \times D } \diff t \otimes Q^{(i)}(\diff r) \otimes \diff s
    < 
    \infty,
    &
    &\forall t_0 \in [0,T].
\end{align*}

The following assumptions will be made throughout, even if they are not explicitly mentioned. 
\begin{enumerate}

    \item[(A1)]
         Let $(\Omega^0,d^0)$ and $(\Omega^1,d^1)$ be complete and separable metric spaces,  endowed with  Borel $\sigma$-algebras  $\Ff^0$ and $\Ff^1$, respectively. We assume that $(\Omega^0, \Ff^0, \PP^0)$ is a complete probability space supporting a $l$-dimensional Poisson process $N$. Similarly, $(\Omega^1, \Ff^1, \PP^1)$ is a complete probability space supporting a $k$-dimensional Wiener process $W$.
				
    \item[(A2)]
		There exists a measurable mapping 
		$X_0:(\Omega^0,\Ff^0) \to (\RR^d,\Bb(\RR^d))$ with law 
		$\nu_0 \in \Pp_2(\RR^d)$; i.e.
		\begin{align*}
            \Esp[^0]{ \Abs{ X_0 }^2 }=\int_{ \RR^d } \Abs{ x }^2 \nu_0( \diff x )
            <
            \infty.
		\end{align*}

\end{enumerate}

\begin{theorem}
\label{Thm:Existence-of-admissible-settings}
Assume that $\mu$ is a random environment on $(\Omega^0,\Ff^0,\PP^0)$, and define $\FF^{\mu,N}:=\FF^{0,\mu}\vee \FF^{0,N}$, where $\FF^{0,N}$ is the filtration generated by
\begin{align*}
    \Ff_t^{0,N} := \sigma \Braces{ 
        N( B \times C \times D ) : \forall (B,C,D) \in \Bb([0,t]) \otimes \Rr \otimes \Bb(\RR_+) 
    }.
\end{align*}
Let $\lambda : \Omega \times [0,T] \times \mathbf{R} \to \RR^l$ satisfy the following:
\begin{enumerate}

    \item
        $\lambda^{(i)}(\cdot,t,r)$ is independent of $\lambda^{(j)}(\cdot,t,r)$ with respect to $\PP^0$ for all $i \neq j$ and all $(t,r) \in [0,T] \times \mathbf{R}$.
    \item 
        $\lambda$ is $\FF^{\mu,N}$-predictable.

    \item
        There exists an increasing sequence of measurable sets $R_j \uparrow R$ such that for all $0 \leq a < b \leq T$,
        \begin{align*}
            &\int_a^b \int_{R_j} \lambda^{(i)}(t,r) Q^{(i)}(\diff r) \diff t < \infty
            &
            &\PP^0-\mbox{a.s.}
        \end{align*}
        for all $R_j$ and for all $i =1, \ldots,l$.
        
\end{enumerate}
Then, there exists an admissible set-up $(\Omega,\Ff,\FF,\PP,X_0,\mu,\lambda,(N^\lambda,\xi),W)$ such that
\begin{itemize}

    \item[(i)]
		$(N^\lambda,\xi)$ is an $\FF^{\mu,N}$-admissible noise on $(\Omega^0,\Ff^0,\PP^0)$  with $\FF^{\mu,N}$-intensity kernel
        \begin{align}
            \label{Eq:K=lambda.Q}
            K^\lambda(t,\diff r) 
            := 
            \lambda(t, r) Q(\diff r) 
            = 
            \Parentheses{ \lambda^{(i)}(t,r) Q^{(i)}(\diff r) }_{i=1, \ldots, l}.
        \end{align}

    \item[(ii)]
		The triplet $(\Ff, \FF, \PP)$ is defined on the product space 
		$\Omega = \Omega^0 \times \Omega^1$.
	
    \item[(iii)]
		$(X_0,\mu,\lambda,(N^\lambda,\xi), W)$  are defined on $\Omega$ as their \emph{natural extensions}:
		\begin{align}
            \label{Eq:(Intro)Natural-Extensions}
            &(X_0,\mu,\lambda,(N^\lambda,\xi))( \omega^0, \omega^1 ) 
                := (X_0,\mu,\lambda,(N^\lambda,\xi))( \omega^0 ),
		  &
            &W( \omega^0, \omega^1 ) := W( \omega^1 ),
		\end{align}
		for all $( \omega^0, \omega^1 ) = \omega \in \Omega$.
	
\end{itemize}
\end{theorem}

\begin{proof}
Let $\Omega := \Omega^0 \times \Omega^1$ and denote by $\Nn$ the set of $\PP^0 \otimes \PP^1$-null sets on $(\Omega,\Ff)$. Define the \emph{space with admissible Poissonian-type noise} as
\begin{align}
	\label{Eq:(Omega,Ff,PP)}
	(\Omega, \Ff, \PP)
	:=
	\Parentheses{
		\Omega^0 \times \Omega^1, 
		(\Ff^{\mu,N}_T \otimes \Ff^{W}_T) \vee \Nn, 
		\PP^0 \otimes \PP^1
	}.
\end{align}
Similarly,  previous filtrations are extended to the product space:
\begin{align}
	\label{Eq:FF^0}
	\FF^0 
	:&= 
	\Parentheses{ \FF^{\mu,N} \otimes \Braces{ \emptyset, \Omega^1 } } \vee \Nn
	=
	\Braces{ 
		( \Ff_t^{\mu,N} \otimes \Braces{ \emptyset, \Omega^1 } ) \vee \Nn=:\Ff^0_t, 
		0 \leq t \leq T 
	},
	\\ \label{Eq:FF^1}
	\FF^1 
	:&= 
	\Parentheses{ \Braces{ \emptyset, \Omega^0 } \otimes \FF^{W} } \vee \Nn
	=
	\Braces{ 
		( \Braces{ \emptyset, \Omega^0 } \otimes \Ff_t^{W} ) \vee \Nn=:\Ff^1_t, 
		0 \leq t \leq T 
	},
	\\ \label{Eq:FF}
	\FF :&= \FF^{0} \vee \FF^{1}.
\end{align}
Note that $\FF$, $\FF^0$ and $\FF^1$ satisfy the usual conditions, and $\FF^0$ and $\FF^1$ are independent. Furthermore, $\FF = \Parentheses{ \FF^{\mu,N} \otimes \FF^{W} } \vee \Nn$.

For the stochastic processes defined on each  coordinates, we  abuse slightly 
of notation and use the same symbols as their natural extension, i.e.
\begin{align}
	\label{Eq:Theta-extended}
	\Theta_t (\omega^0, \omega^1) 
	&:= 
	\Theta_{t} (\omega^0)\ 
	&
	&\forall (t, \omega^1) \in [0,T] \times \Omega^1,
	\\ \nonumber
	&\Parentheses{\mbox{resp. }\Theta_{t}(\omega^1)},
	&
	&\Parentheses{\mbox{resp. }\forall (t, \omega^0) \in [0,T] \times \Omega^0}.
\end{align}
Since  $\Omega^i$ are  Polish spaces, with $\Ff^i_T$ the corresponding Borel $\sigma$-algebra, the existence of a regular conditional probability measure guarantees that $\PP$ admits a disintegration \cite{bogachev_measure_2007} of the form 
\begin{align}
	\label{Eq:desintegration-of-PP}
	\PP\Brackets{ \diff \omega^0, \diff \omega^1 }
	=
	\PP^0\Brackets{ \diff \omega^0 } \otimes P^1 \Parentheses{ \omega^0, \diff \omega^1 },
\end{align}
for some probability measure $\PP^0$ on $\Parentheses{ \Omega^0, \Ff_T^0 }$ and $P^1$  a stochastic kernel 
on $\Parentheses{ \Omega^1, \Ff_T^1 }$ given 
$\Parentheses{ \Omega^0, \Ff_T^0 }$, such that  the following claims are equivalent \cite[Proposition 1.10]{carmona_probabilistic_2018-1}:
\begin{enumerate}

	\item
		For each $t \in [0,T]$ and $A \in \Ff_t^1$, the mapping
		\begin{align*}
			\omega^0 \mapsto P^1(\omega^0, A)
		\end{align*}
		is measurable with respect to  $\Ff_t^0$.
		
	\item
		Each martingale $M$ on 
		$\Parentheses{ \Omega^0, \Ff_T^0, \FF^{\mu,N}, \PP^0 }$
		extended to $\Omega = \Omega^0 \times \Omega^1$ in the same way as 
		\eqref{Eq:Theta-extended} remains a martingale on 
		$\Parentheses{ \Omega, \Ff, \FF, \PP }$.
\end{enumerate}

For the Admissible Common Noise, it follows directly from \cite[Thm. 5.7.3]{bremaud_point_2020} that there exists an $\FF^{\mu,N}$-adapted marked Poisson process $(N^\lambda,\xi)$ with $\FF^{\mu,N}$-intensity kernel $K^\lambda$ as in \eqref{Eq:K=lambda.Q}, characterized by the lifting
\begin{align*}
    \eta^\lambda(\diff t, \diff r)
    :=
    N( \diff t, \diff r, [0,\lambda(t,r)] ).
\end{align*}
For every $C \in \Rr$, let $\tilde{\eta}^{\lambda,C}$ be the compensated martingale of $\eta^{\lambda,C}=\eta^\lambda([0,\cdot],C)$:
\begin{align}
	\label{Eq:tilde-eta^lambda}
	\tilde{\eta}^{\lambda,C}_t
	:&=
	\eta^{\lambda,C}_t - \int_{[0,t] \times C} K^\lambda(s,\diff r) \diff s,
	&
	& 0 \leq t \leq T.
\end{align}
Then,  $\tilde{\eta}^{\lambda,C}$ and $W$ are martingales on $(\Omega,\Ff,\PP)$ with respect to $\FF^0$ and $\FF^1$, respectively, and therefore we get the independence and joint compatibility with respect to $\FF$. Since $\tilde{\eta}^{\lambda,C}$ is an $\FF$-martingale,  the compensation formula holds; see \cite{bremaud_point_2020}. 

For the Idiosyncratic Noise, recall that the quadratic variation of a process does not depend on the filtration. In other words, $W$ is a local $\FF$-martingale such that $W_0 = 0$ and $\Brackets{W,W}_t = t$. By Lévy's characterization of the Brownian motion, the result follows.

\end{proof}

The admissible set-up $(\Omega,\Ff,\FF,\PP,X_0,\mu,\lambda,(N^\lambda,\xi),W)$ introduced in the previous theorem remains fixed throughout. In order to state the existence and uniqueness of solutions to the system \eqref{Eq:Forward-Component}-\eqref{Eq:Initial-Terminal-Conditions}, the following standing assumptions are needed.\\[.1cm] 

\noindent
{\bf Assumption A}
\begin{enumerate}
    \item
		For all $
			(\nu, x, y, z, u )
			\in
			\Pp_2(\RR^d) 
				\times \RR^d
				\times \RR^n 
				\times \RR^{n \times k} 
				\times \mathfrak{H}^\lambda
		$, the process
		\begin{align*}
            (b, f, \sigma)( \cdot, \nu, x, y, z, u ) 
		\end{align*}
		take values in $\RR^d \times \RR^n \times \RR^{d \times k}$, is $\FF$-progressively measurable, and
        \begin{align*}
            (b, f, \sigma)(\cdot,\mu_\cdot,0,0,0,0)
		  \in
            \HH^2(W,\tilde{\eta}^\lambda) \times \HH^2(W,\tilde{\eta}^\lambda) \times \HH^2(W) 
        \end{align*}

    \item
        For each $
		  (\nu, x, y, z, u )
		  \in
            \Pp_2(\RR^d) \times \RR^d \times \RR^n  \times \RR^{n \times k} \times \mathfrak{H}^\lambda,
		$ the process
		\begin{align*}
		  \gamma( \cdot, \nu, x, y, z, u )
		\end{align*}
        takes values in $\RR^{d \times l}$, is $\FF \otimes \Rr$-predictable, and
        \begin{align*}
            \gamma(\cdot,\mu_\cdot,0,0,0,0) \in \HH^2(\tilde{\eta}^\lambda);
        \end{align*}

    \item 
        $g(\nu,x) \in L^2(\Ff_T)$ for all $(\nu,x) \in \Pp_2(\RR^d) \times \RR^d$;
    
    \item
		There exist $\FF$-predictable processes
		\begin{align}
		  \label{Eq:(Intro)Lipschitz-Coefficients}
		  K, \underline{K}, K^W, K^\lambda 
		  \ &:\  
		   \Omega \times \RR_+ 
		  \to 
		   \RR_+,
		\end{align}
        and a positive constant $K^0 > 0$, such that for any $x,x' \in \RR^d$, $\nu,\nu' \in \Pp_2(\RR^d)$, $y,y' \in \RR^d$, $z,z' \in \RR^{ d \times m }$ and $u,u' \in \mathfrak{H}^\lambda$, the following inequality holds for $\diff \PP \otimes \diff t$-almost all $(\omega,t) \in \Omega \times [0,T]$:
		\begin{align*}
            &\Abs{ 
				\delta (b,f)_t( \nu, x, y, z, u_t(\omega) )(\omega)
		  }^2
            +
            \Norm{ 
				\delta \sigma_t ( \nu, x, y, z, u_t(\omega) )(\omega)
		  }^2_{\mathrm{Fr}}
            \\
            &\qquad+
            \Norm{
				\delta \gamma_t ( \cdot, \nu, x, y, z, u_t(\omega) )(\omega)
		  }^2_{\lambda_t(\omega) }
            + 
            \Abs{ \delta g(\nu,x)(\omega) }^2
		  \\
		  &\leq
		  K^0\ W_2(\nu,\nu')^2
		  +
		  K_t(\omega) \Abs{ \delta x }^2
		  +
		  \underline{K}_t(\omega) \Abs{ \delta y }^2
		  +
		  K^W_t(\omega) \Norm{ \delta z }_{ \mathrm{Fr} }^2
		  \\
		  &\qquad
		  +
		  K^\lambda_t(\omega) \Norm{ 
				\delta u_t(\omega,\cdot)
		  }_{ \lambda_t ( \omega ) }^2.
		\end{align*}

    \item
        Furthermore, the Lipschitz coefficients $K, \underline{K}, K^W, K^\lambda$
        from \eqref{Eq:(Intro)Lipschitz-Coefficients} are 
        $\diff \PP \otimes \diff t$-a.e. bounded from below away from zero by a positive constant 
        $K_*>0$, and from above by another positive constant $K^* $.
		
    \item 
        There exists a full rank $\RR^{n \times d}$ matrix $G$ such that, for $(\omega,t)\in \Omega \times [0,T]$ fixed, the mapping $A^\mu: \RR^d \times \RR^n \times \RR^{d \times k}\times \mathfrak{H}^\lambda\to \RR^d \times \RR^n \times \RR^{d \times k} \times \mathfrak{H}^\lambda$, defined as
        \begin{align}
        \label{DefA}
            A_t^\mu(x,y,z,u)(\omega)
            :&=
            \begin{bmatrix}
                -G^\top f_t^\mu(x,y,z,u_t)(\omega)
                \\
                G b_t^\mu(x,y,z,u_t)(\omega)
                \\
                G \sigma_t^\mu(x,y,z,u_t)(\omega)
                \\
                G \gamma_t^\mu(\cdot,x,y,z,u_t)(\omega)
            \end{bmatrix},
        \end{align}
        exists and satisfies the $G$-{\it monotonicity condition} described next. Namely, there exist nonnegative constants $\beta_1$, $\beta_2$, $\beta_3$ with $\beta_1 + \beta_2 > 0$ and $\beta_2 + \beta_3 > 0$ such that the following inequalities hold $\diff \PP \otimes \diff t$-a.e.:
        \begin{align*}
            \Big( \delta {A}_t^\mu(x,y,z,u_t) \cdot \delta (x, y, z, u_t) \Big) 
            \leq
            &- 
            \beta_1 \Abs{ G ( \delta x ) }^2
            \\
            &-
            \beta_2 \Braces{ 
                \Abs{ G^\top ( \delta y ) }^2 + \Norm{ G^\top ( \delta z ) }^2_{\mathrm{Fr}} + \Norm{ G^\top ( \delta u_t ) }^2_{\lambda_t}
            },
        \end{align*}
        where inner product at the left hand side corresponds to the sum of inner products in the corresponding (matrix) spaces, and
        \begin{align*}
            &\Big( \delta g_T^\mu (x)(\omega) \cdot G ( \delta x ) \Big)
            \geq 
            \beta_3 \Abs{ G ( \delta x ) }^2,
        \end{align*}
        for all $ (x,y,z,u), (x',y',z',u') \in \RR^d \times \RR^n \times \RR^{ d \times k } \times \mathfrak{H}^\lambda$.

\end{enumerate}

\begin{remark}
\label{dimen}
Depending on the dimension of $G$, we have that either $\beta_1,\beta_3 > 0$ ($n > d$), or $\beta_2>0$ ($n \leq d$). Indeed, take $\delta(x,y,z,u) \neq 0$ and assume that the coefficients are non-constant; otherwise, the assumption holds trivially. If $n \geq d$, $G$ is full column rank and $G (\delta x) \neq 0$, and hence, from the definition of $G$-monotonicity, $\beta_1 = 0$ (respectively $\beta_3 = 0$) if and only if $-A$ (resp. $g$) is a constant function. An analogue result can be obtained for $\beta_2$ when $n \leq d$. Moreover, when 
$n=d$, we can always take $G = \mathrm{Id}_n$.
\end{remark}

The main result of this work is presented below, postponing its proof to the next sections.
\begin{theorem}
\label{Thm:Existence-Uniqueness-FBSDE}
There exists a unique process 
\begin{align*}
    ( X,Y,Z,U,M)
    \in \Hh^2(\RR^d)
        \times \Hh^2(\RR^n)
        \times \HH^{2}(W)
        \times \HH^{2}(\tilde{\eta}^\lambda) \times
        \Parentheses{ \Ll^2(W) \oplus \Ll^2(\tilde{\eta}^\lambda)} ^{\perp}
\end{align*}
(up to $\PP$-indistinguishability)  satisfying the FBSDE
\begin{align}
    \tag{\ref{Eq:Forward-Component}}
    X_t^\mu
    =&
    X_0^\mu
    +
    \int_0^t b^\mu_s \Parentheses{ X^\mu, Y^\mu, Z^\mu, U^\mu } \diff t
    +
    \int_0^t \sigma^\mu_s \Parentheses{ X^\mu, Y^\mu, Z^\mu, U^\mu } \diff W_t
    \\ \nonumber
    &+
    \int_{ [0,t] \times \mathbf{R} } \gamma^\mu_s \Parentheses{ r,
        X^\mu, Y^\mu, Z^\mu, U^\mu, r
    } \tilde{\eta}^\lambda( \diff s, \diff r ),
    \\ \tag{\ref{Eq:Backward-Component}}
    Y_t^\mu
    =&
    Y_T^\mu
    +
    \int_t^T f^\mu_s \Parentheses{ 
        X^\mu, Y^\mu, Z^\mu, U^\mu
    } \diff s
    -
    \int_t^T Z^\mu_s \diff W_s
    \\ \nonumber
    &-
    \int_{(t,T] \times \mathbf{R}} U_{s}^\mu(r) \tilde{\eta}^\lambda( \diff s, \diff r )
    -
    (M^\mu_T - M^\mu_t),
    \\ \tag{\ref{Eq:Initial-Terminal-Conditions}}
    X_0^\mu \sim& \mu_0 \in \Pp_2(\RR^d; \Ff_0),
    \qquad
    Y_T^\mu = g_T^\mu(X) \in L^2(\RR^n; \Ff_T).
\end{align}
\end{theorem}

The choice of  the $G$-monotonicity condition is based on ideas adapted from equilibrium problems; see \cite{oettli_existence_1998}, \cite{brokate_generalized_nodate}. These problems arise naturally in optimization problems, Nash equilibria of noncooperative games and variational inequalities. In the context of FBSDEs, this condition has been also useful to study fully coupled systems of stochastic differential equations with both initial and terminal conditions; see, for instance,  \cite{peng_fully_1999} \cite{zhen_forward-backward_1999}.   The implementation of the \emph{method of continuation} for solving the system of FBSDEs 
studied in this work depends heavily on this condition.

\subsection{Hamiltonian systems}

The following Hamiltonian system is inspired by stochastic control problems. We present the unmarked Poisson process in one dimension for ease of notation. Given an admissible set-up, let $b,\sigma,f$, $\gamma$ and $g$ be the elements of the system, with
\begin{align*}
    b,\sigma,f,\gamma &: \Omega \times [0,T] \times \Pp_2(\RR) \times \RR \to \RR,
    \\ 
    g&:\RR \times \Pp_2(\RR) \to \RR,
\end{align*}
satisfying suitable technical assumptions. To avoid degeneracies on the jumps, we assume $0 \not\in \mathrm{Image}(\gamma)$. 

Define the Hamiltonian $\Omega \times [0,T] \ni (\omega,t) \mapsto H(\omega,t; \nu, x, y, z, u)$ by
\begin{align*}
    H(\omega, t, \nu, x, y, z, u)
    :=&
    f(\omega, t, \nu, x)
    +
    y b(\omega, t, \nu, x)
    +
    z \sigma (\omega, t, \nu, x)
    +
    u \lambda_t(\omega) \gamma(\omega, t, \nu, x),
\end{align*}
for  $ (\nu,x,y,z,u) \in \Pp_2(\RR) \times \RR^4$. Let us write by abbreviation 
\begin{align*}
	H_t^\mu(\omega) 
	= 
	H( \omega, t; 
		\mu_{t-}(\omega), 
		X_t(\omega), 
		Y_t(\omega), 
		Z_t(\omega), 
		U_t(\omega) 
	).
\end{align*}
The Jacobian of $H$, assuming that the derivatives involved exist, has the form
\begin{align*}
    \Diff H_t^\mu &( x, y, z, u )
    \\
    &=
    \begin{bmatrix}
		\partial_x f_t^\mu(x) + y \partial_x b_t^\mu(x) + z \partial_x \sigma_t^\mu(x) 
            + u \lambda_t \partial_x \gamma_t^\mu(x) 
		&
		b_t^\mu(x)
		&
		\sigma_t^\mu(x)
        &
        \lambda_t \gamma_t^\mu(x)
    \end{bmatrix}.
\end{align*}

The maximum principle in stochastic control is composed by a forward system, representing the evolution of a controlled system, and a backward equation, corresponding to the {\it adjoint equation} \cite{oksendal_applied_2019}.  
\begin{align*}
    X_t
    =&
    X_0
    +
    \int_0^t \partial_y H^\mu_s \diff s
    +
    \int_0^t \partial_z H^\mu_s \diff W_s
    +
    \int_0^t \frac{1}{\lambda_t} \partial_u H^\mu_s \diff \tilde{N}^\lambda s,
    \\
    Y_t 
    =&
    \partial_x g_T^\mu(X)
    +
    \int_t^T \partial_x H^\mu_s \diff s
    -
    \int_t^T Z_s \diff W_s
    -
    \int_{t+}^T U_{s} \diff \tilde{N}^\lambda s - (M_T - M_t).
\end{align*}
Thus, within the context of an exogenous  random environment, the existence and uniqueness of solution to this FBSDE   had to be solved in order to applied the maximum principle approach to solve the optimal control problem behind it. In this case,  the corresponding process $A^\mu$ is given by the gradient of $H$, i.e. $A^\mu = \Diff H^{\mu\top}$.

\section{Preliminaries from uncoupled forward and backward equations }
\label{PUFBE}
 When the forward and backward components of the system of interest are decoupled, their existence and uniqueness can be obtained  from the existing literature, analyzing independently each component, under a set of simplified but consistent assumptions.  

More precisely, define the \emph{SDE in the (\textit{fixed}) random environment $\mu$, with coefficients $(b, \sigma, \gamma)$ and initial condition $X_0$}, as the solution of the integral equation 
\begin{align}
    \label{Eq:(FSDE)X=b.dt+sigma.dW+gamma.deta}
    X_t
    =&
    X_0
    +
    \int_0^t b(s, \mu_s, X_{s} ) \diff s
    +
    \int_0^t \sigma(s, \mu_s, X_{s} ) \diff W_s
    +
    \int_{[0,t] \times \mathbf{R}} \gamma(s, r, \mu_{s-}, X_{s-} ) 
        \tilde{\eta}^\lambda (\diff s, \diff r),
\end{align}
for all $t \in [0,T]$, where 
\begin{align}
    \label{Eq:General-b-sigma}
    &(b, \sigma)
    \ :\ 
    \Omega \times [0,T] \times \Pp_2(\RR^d) \times \RR^d
    \to
    \RR^d \times \RR^{d \times k},
    \\
    &\gamma
    \ :\ 
    \Omega \times [0,T] \times \mathbf{R} \times \Pp_2(\RR^d) \times \RR^d
    \to
    \RR^{d \times l}. \nonumber
\end{align}
Here, the r.h.s. of (\ref{Eq:(FSDE)X=b.dt+sigma.dW+gamma.deta}) should be understood as a component-wise $\PP$-a.s. equality
\begin{align*}
    &X_0^{(i)}
    +
    \int_0^t b_s^{(i,1)} \diff s
    +
    \int_0^t \sigma_s^{(i, \cdot )} \diff W_s
    +
    \int_{[0,t]\times \mathbf{R}} \gamma_{s}^{(i, \cdot )}(r) \tilde{\eta}^\lambda( \diff s, \diff r )
\end{align*}
for each $i = 1, \ldots, d$, where the sup-fix $(i,\cdot)$ denotes the $i$-th row of the corresponding matrix. To ease the notation,  we  write the dependence on the time and environment as
\begin{align*}
    (b,\sigma)_t^\mu(X)(\omega)
    :&=
    (b,\sigma)(\omega, t, \mu_t(\omega), X_t(\omega) ),
    &
    \gamma_t^\mu(r,X)(\omega)
    :&=
    \gamma(\omega, t, r, \mu_{t-}(\omega), X_{t-}(\omega) ),
\end{align*}
for all $(\omega,t) \in \Omega \times [0,T]$. The corresponding simplified versions of Assumption A  (i.e. by taking a null backward component, $G$-monotonicity implies linear growth), allow us to  apply  \cite[Theorems 14.21-14.23]{jacod_calcul_2006}, as well as the usual martingale results and Gr\"onwall's inequality, in order to obtain the following result. 

\begin{proposition}
Let $(b,\sigma),\gamma$ be defined as in \eqref{Eq:General-b-sigma}. Under Assumption A and taking $(f,g)$ identically zero, there exists a unique $\FF$-adapted process $X$, defined on $(\Omega, \Ff, \PP)$ with   $\Norm{ X }_{*} < \infty$,  solving the SDE  \eqref{Eq:(FSDE)X=b.dt+sigma.dW+gamma.deta}.
\end{proposition}

Solving the backward component in the presence of a random environment requires taking into account the additional information contained in $\FF$, since its solution cannot depend on the future information of process $\mu$. However, the nature of the noise we are considering admits the possibility of applying the existing theory of BSDEs.

Let $\tilde{M} \in \Hh^2(\RR^n)$   be an $\FF$-adapted martingale defined on  $(\Omega, \Ff, \PP)$,  admiting a decomposition of the form $\tilde{M} = \tilde{M}_0 + \tilde{M}^c +\tilde{M}^d$, with $\tilde{M}_0$  an $\Ff_0$-measurable random variable, $\tilde{M}^c$  a square-integrable martingale with continuous paths and $\tilde{M}^d$  a purely discontinuous, square integrable martingale.  Intuitively $\tilde{M}^c$ and $\tilde{M}^d$ are related to $W$ and $(N^\lambda,\xi)$, respectively, by martingale representation theorems. However, in our case, these processes do not represent all the randomness affecting the system, and therefore, in order to ensure that martingale $\tilde{M}$ is $\Ff_t$-measurable at each time $t$, an additional term must be added.

\begin{remark}
In order to address the previous claim, we now present a couple of known results from the theory of BSDEs; see  \cite[Propositions 2.5, 2.6, and Corollary 2.7]{papapantoleon_existence_2018}. The results shown there are presented in a general setting, however, for the sake of completeness and readability, we reproduce them here in the form of a single proposition, adapted to our framework and notation, restricting ourselves to  the case of marked Poisson point processes.
\end{remark}

\begin{proposition} 
\label{Prop:Orthogonal-Decomposition-in-Admissible-Set-up}
Each random variable $\zeta \in L^2(\Ff_T; \RR^n)$ has an orthogonal decomposition with respect to $(N^\lambda,\xi)$ and $W$; namely, there exists a unique (up to $\PP$-indistinguishability) triplet $(Z,U,M)$ such that:
\begin{enumerate}
    
    \item
		$Z \in \HH^2(W)$ and $U \in \HH^2(\tilde{\eta}^\lambda)$;
	
    \item
		the stochastic processes:
        \begin{align*}
            &\Parentheses{ \int_0^t Z_s \diff W_s, t \in [0,T] },
            &
            &\Parentheses{ \int_{[0,t] \times \mathbf{R} } U_s(r) \tilde{\eta}^\lambda(\diff s, \diff r), t \in [0,T] };
        \end{align*}
		are   $\FF$-adapted and independent  martingales;
    \item
		$M \in \Hh^2(\RR^d)$, with $M_0=0$, and orthogonal to $(W, \tilde{\eta}^\lambda)$;
	
    \item 
        the random variable $\zeta$ can be written as follows:
        \begin{align*}
            \zeta
            =
            \Esp{ \zeta \middle| \Ff_0 }
            +
            \int_0^T Z_s \diff W_s
            +
            \int_{[0,T] \times \mathbf{R} } U_s(r) \tilde{\eta}^\lambda(\diff s, \diff r)
            +
            M_T.
	\end{align*}
	
\end{enumerate}

\end{proposition}

Proposition \ref{Prop:Orthogonal-Decomposition-in-Admissible-Set-up} provides the first elements to write BSDEs with a terminal condition $\zeta \in L^2(\Ff_T; \RR^n)$, defined on an admissible set-up, in terms of  an additional martingale $M$. Namely,
\begin{align}
    \label{Eq:(BSDE)Y-complete}
    Y_t
    =&
    \zeta
    +
    \int_t^T f_s \diff s
    -
    \int_t^T Z_s \diff W_s
    -
    \int_{(t,T] \times \mathbf{R} } U_s(r) \tilde{\eta}^\lambda ( \diff s, \diff r )
    - 
    (M_T - M_t),
\end{align}
where $f$ is the driver of the system, an $\FF$-predictable process. When $f$ has the form
\begin{align*}
    f_s(\omega) = f( \omega, s, \mu_s(\omega), Y_s(\omega), Z_s(\omega), U(\omega, s, \cdot) ),
\end{align*}
we say that \eqref{Eq:(BSDE)Y-complete} is a \emph{BSDE on the random environment $\mu$}.

\begin{remark}
Just as in the forward case, \eqref{Eq:(BSDE)Y-complete} should be understood as a handy matrix notation for the component-wise identity
\begin{align*}
    Y_t^{(i)}
    =&
    \zeta^{(i)}
    +
    \int_t^T f^{(i,1)}_s \diff s
    -
    \int_t^T Z_s^{(i, \cdot )} \diff W_s
    -
    \int_{(t,T] \times \mathbf{R} } U^{(i,\cdot)}_s(r) \tilde{\eta}^{\lambda} ( \diff s, \diff r^j )
    -
    (M_T^{(i)} - M_t^{(i)}).
\end{align*}
\end{remark}

As before,  this particular form of driver $f$, and the choice made for noise process, allows us to rely on existent results of the general theory of  BSDEs. Since the predictable quadratic variation of both $W$ and $\tilde{\eta}^\lambda$ is absolutely continuous with respect to the Lebesgue measure on $[0,T]$, and an orthogonal martingale decomposition holds, the existence of a unique solution process $(Y,Z,U,M)$ can be obtained  from \cite[Theorem 3.5 and Corollary 3.6]{papapantoleon_existence_2018}, whenever  $f$ satisfies 
the following conditions:
\begin{enumerate}
    \item 
        \textbf{Measurability} For every $
            (\nu, y, z, u) 
            \in 
            \Pp_2(\RR^d) \times \RR^n \times \RR^{ n \times k } \times \mathfrak{H}^\lambda$ the mapping 
            $(\omega,t) \mapsto f_t( \omega, \nu, y, z, u )$ is $\FF$-progresively measurable;
		
    \item 
        \textbf{Lipschitz} There exist a triplet $(\underline{K}, K^W, K^\eta)$ of positive, $\FF$-predictable processes and a  non-negative constant $K^0 \geq 0$ such that, for any $\nu,\nu' \in \Pp_2(\RR^d)$, $y,y' \in \RR^d$, $z,z' \in \RR^{ d \times m }$ and $u,u' \in \mathfrak{H}^\lambda$, the following inequality holds $\diff \PP \otimes \diff t$-a.e. $(\omega,t) \in \Omega \times [0,T]$:
		\begin{align*}
            &\Abs{ \delta f_t( \nu, y, z, u )(\omega) }^2
            \\
            &\quad\leq
            K^0\ W_2(\nu,\nu')^2
            +
            \underline{K}_t(\omega) \Abs{ \delta y }^2
            +
            K^W_t(\omega) \Norm{ \delta z }_{ \mathrm{Fr} }^2
            +
            K^\lambda_t(\omega) \Norm{ \delta u_t(\omega, \cdot ) }_{ \lambda_t ( \omega ) }^2.
		\end{align*}
\end{enumerate}
Moreover,  for some sufficiently large $\beta^* > 0$, 
\begin{align*}
    &\Esp{ e^{\beta^*A_T} \Abs{\zeta}^2 } < \infty
    &
    &\mbox{and}
    &
    &\Esp{ \int_0^T e^{\beta^*A_t} \frac{\Abs{f_t^\mu(0,0,0)}}{ \alpha^2_t } \diff t } < \infty,
\end{align*}
where $A$ and $\alpha$ are non-negative processes defined as
\begin{align}
    \label{Eq:alpha-A}
    \alpha_t^2 :&= \max\Braces{ \sqrt{\underline{K}_t}, K^W_t, K^{\eta}_t },
    &
    A_t :&= \int_0^t \alpha^2_s \diff s.
\end{align}
\begin{proposition}[\cite{papapantoleon_existence_2018}]\label{PropIn}
Let $
    f : \Omega \times [0,T] \times \Pp_2(\RR^d) \times \RR^n \times \RR^{ n \times k } 
        \times \mathfrak{H}^\lambda
    \to
    \RR^n
 $ be a driver satisfying the above Measurability and Lipschitz conditions. Then, there exists a unique (up to $\PP$-indistinguishability) process $(Y,Z,U,M)$ with
 \begin{align*}
     &Y,M : \Omega \times [0,T] \to \RR^n,
     &
     &Z : \Omega \times [0,T] \to \RR^{ n \times k },
     &
     &U : \Omega \times [0,T] \times \mathbf{R} \to \RR^{ n \times l }
 \end{align*}
 such that 
\begin{itemize}
    \item 
        $Y$ and $M$ are $\FF$-adapted, $Z$ is $\FF$-progressively measurable, and $U$ is $\FF \otimes \Rr$-predictable;

    \item 
        $(Y,Z,U,M)$ is a solution process to the BSDE in random environment \eqref{Eq:(BSDE)Y-complete};

    \item 
        $\Norm{ (Y,Z,U,M) }_{*,\beta^*,A}^2 < \infty$, with
        \begin{align*}
            \Norm{ (Y,Z,U,M) }_{*,\beta^*,A}^2
            :=
            \EE\Bigg[
                &\sup_{t\in[0,T]} \Braces{ e^{\beta^*A_t} \Abs{Y_t}^2 } 
                +
                \int_0^T e^{\beta^*A_t} \Braces{ 
                    \Norm{Z_t}^2_\mathrm{Fr} 
                    + 
                    \Norm{ U_t(\cdot) }_{\lambda_t}^2 }
                \diff t
                \\
                &+
                \int_0^T e^{\beta^*A_t} \diff \mathrm{tr} \Angles{M}_t
            \Bigg].
        \end{align*}
\end{itemize}

\end{proposition}
The conditions adopted in the previous result are perfectly compatible with our standing  Assumption A, as it is stated  next.    

\begin{corollary}
\label{Cor:Bounded-Lipschitz-Driver}Under Assumption A and by taking  $(b,\sigma,\gamma)$ identically zero with $g\equiv \zeta\in L^2(\Ff_T; \RR^n)$, the BSDE \eqref{Eq:(BSDE)Y-complete} has a unique $\FF$-adapted solution $(Y,Z,U,M)$ in $\Ss^2 \times \HH^{2}(W) \times \HH^{2}(\tilde{\eta}^\lambda) \times \Hh^{2,\perp}$. Furthermore, the norms $\Norm{\ \cdot\ }_{ \ast }$ and $\Norm{\ \cdot\ }_{ \ast, \beta, A }$ are equivalent for each $\beta \geq 0$.
\end{corollary}

\begin{proof}
Observe that 
\begin{align*}
	&e^{\beta K_*T} \Norm{ \zeta }_{ L^2 }
	\leq
	\Norm{ \zeta }_{ L^2_{\beta, A} }
	\leq
	e^{\beta K^*T} \Norm{ \zeta }_{ L^2 }
	<
	\infty,
\end{align*}
and
\begin{align*}
    \frac{ e^{\beta K_*T} }{K^*} \Esp{ \int_0^T \Abs{f_t^\mu(0,0,0)}\diff t }
    &\leq
    \Esp{ \int_0^T e^{\beta A_t} \frac{\Abs{f_t^\mu(0,0,0)}}{ \alpha^2_t } \diff t }
    \\
    &\leq
    \frac{ e^{\beta K^*T} }{ K_* } \Esp{ \int_0^T \Abs{f_t^\mu(0,0,0)}\diff t } 
    <
    \infty
\end{align*}
$\diff \PP \otimes \diff t$-a.e. $(\omega,t)\in \Omega \times [0,T]$, for all $\beta \geq 0$. Thus, for any $\beta > 0$, $((W,\tilde{\eta}^\lambda), \FF, T, \zeta, f, \diff t)$ is a standard data under $\beta$ in the sense of    \cite[Definition 3.2]{papapantoleon_existence_2018}. Taking $\beta^*$ sufficiently large, there exists a unique $
    (Y,Z,U,M) 
    \in
    \Ss^2_{\beta^*,A} \times \HH^{2,W}_{\beta^*,A} \times \HH^{2,N^\lambda}_{\beta^*,A} 
    \times \Hh^{2,\perp}_{\beta^*,A}
$  satisfying \eqref{Eq:(BSDE)Y-complete}. 

The equivalence between norms is due to the estimate 
\begin{align*}
    &0 < K_* t \leq A_t(\omega) \leq K^* T,
    &
    &\diff \PP \otimes \diff t-\mbox{a.e.}
\end{align*}
Indeed, for each $\beta > 0$, we have that
\begin{align*}
    \Norm{ (Y,Z,U,M) }^2_{ \ast }
    &=
    \Norm{ Y }_{ \Ss^2 }^2 + \Norm{ Z }_{ \HH^{2}(W) }^2 + \Norm{ U }_{ \HH^{2}(N^\lambda) }^2 
        + \Norm{ M }_{ \Hh^{2,\perp} }^2
    \\
    &\leq
    \Norm{ Y }_{ \Ss^2_{ \beta, A } }^2 + \Norm{ Z }_{ \HH^{2,W}_{ \beta, A } }^2 
        + \Norm{ U }_{ \HH^{2,N^\lambda}_{ \beta, A } }^2 
        + \Norm{ M }_{ \Hh^{2,\perp}_{ \beta, A }  }^2
    \\
    &=
    \Norm{ (Y,Z,U,M) }^2_{ \ast, \beta, A } 
    < 
    \infty.
\end{align*}
Moreover,
\begin{align*}
    &\Norm{ (Y,Z,U,M) }^2_{ \ast, \beta, A } 
    =
    \Norm{ Y }_{ \Ss^2_{ \beta, A } }^2 + \Norm{ Z }_{ \HH^{2,W}_{ \beta, A } }^2 
        + \Norm{ U }_{ \HH^{2,\tilde{\eta}^\lambda}_{ \beta, A } }^2 
        + \Norm{ M }_{ \Hh^{2,\perp}_{ \beta, A }  }^2
    \\
    &\qquad\leq
    \Norm{ e^{\beta K^* T} Y }_{ \Ss^2 }^2 + \Norm{ e^{\beta K^* T} Z }_{ \HH^{2}(W) }^2 
        + \Norm{ e^{\beta K^* T} U }_{ \HH^{2}(\tilde{\eta}^\lambda) }^2 
        + \Norm{ e^{\beta K^* T} M }_{ \Hh^{2,\perp} }^2
    \\
    &\qquad=
    e^{\beta K^* T} \Norm{ (Y,Z,U,M) }^2_{ \ast }
    <
    \infty.
\end{align*}

\end{proof}

The following \emph{a priori} estimate, obtained by applying Lemma 3.4 and Proposition 3.13 of \cite{papapantoleon_existence_2018} to our framework, will be useful in the next section.

\begin{lemma}
\label{Lemma:(BSDE)A-priori-estimates}
Under the same conditions as Corollary \ref{Cor:Bounded-Lipschitz-Driver}, let $(Y,Z,U,M)$ and $(Y'$, $Z'$, $U'$, $M')$ be the corresponding solutions to \eqref{Eq:(BSDE)Y-complete} with driver $f^1$ and $f^2$, and terminal conditions $\zeta$ and $\zeta'$, respectively. Then, there exists $\beta^* > 0$ such that the estimate
\begin{align*}
    \Norm{ \delta(Y,Z,U,M) }_{ \ast }^2
    &\leq 
    c_\ast(f^1,f^2,\beta,T) \Norm{ \delta \zeta }^2_{ L^2 }
    +
    c_{**}(f^1,f^2,\beta,T) \Esp{ \int_0^T \Abs{\delta f_t^\mu(Y',Z',U')}  \diff t },
\end{align*}
holds for all $\beta \geq \beta^*$, where $c_*(f^1,f^2,\beta,T)$ and $ c_{**}(f^1,f^2,\beta,T)$ are positive constants, and $\delta f^{\mu}(Y',Z',U') := f^{1\mu}(Y',Z',U') - f^{2\mu}(Y',Z',U')$.

\end{lemma}

\begin{proof}
Let $K^*,K_* \in \RR_+$ be upper and lower bounds for the Lipschitz coefficients of $\delta f = f^1 - f^2$, respectively. From the proof of Corollary \ref{Cor:Bounded-Lipschitz-Driver}, the Lipschitz contitions hold for $f^1$ and $f^2$ with the same $\beta > 0$, and hence for $\delta f$. 

Fix $\beta > 0$, and let $\alpha^*$ and $A^*$ be defined as in \eqref{Eq:alpha-A} for the driver $\delta f$. Since $\Angles{W}$ and $\Angles{\tilde{\eta}^{\lambda,\mathbf{R}}}$ are absolutely continuous with respect to the Lebesgue measure, which in turn is atomless, there exists a constant $\Phi > 0$ such that $(A^*_t - A^*_{t-}) \leq \Phi$ $\diff t \otimes \diff P$-a.e. Furthermore: in this case the constant $\Phi$ can be chosen independently from the rest of the parameters. 

Let $\mathrm{M}^{\Phi}_*(\beta)$ be the function defined in Lemma 3.4 of \cite{papapantoleon_existence_2018}. By definition,
\begin{align*}
    \lim_{\beta \to \infty} \mathrm{M}^{\Phi}_*(\beta) = 9e\Phi,
\end{align*}
and since $\Phi$ can be chosen independently, for any $\Phi < \frac{1}{2}$ there exists $\beta^*(\Phi)$ such that for all $\beta \geq \beta^*$ the conditions from Proposition 3.13 of \cite{papapantoleon_existence_2018} hold. Without loss of generality, take $\Phi=(36e)^-1$; then, from Proposition 3.13 of \cite{papapantoleon_existence_2018}, there exists $\beta^* > 0$ (determined independently of the parameters) and two positive constants $c_1(\beta)$ and $c_2(\beta)$, depending only on $\beta$, such that for all $\beta > \beta^*$, following estimates hold:
\begin{align*}
    \Norm{ \delta(Y,Z,U,M) }_{ \ast }^2
    \leq &
    \Norm{ \delta(Y,Z,U,M) }_{ \ast, \beta, A^* }^2
    \\
    \leq &
    c_1(\beta) \Norm{ \delta \zeta }^2_{ L^2_{\beta, A^* } }
    +
    c_2(\beta) 
        \Esp{ \int_0^T e^{\beta A^*_t} \frac{\Abs{\delta f_t^\mu(Y',Z',U')}}{ \alpha^{*2}_t } \diff t }.
\end{align*}
Lastly, since $\alpha^*$ and $A^*$ depend on the parameters only through the Lipschitz coefficients of the drivers (not on $\zeta$ and $\zeta'$), the result is obtained by observing that the last term of the previous inequality is bounded by
\begin{align*}
    \underbrace{ c_1(\beta) e^{\beta K^* T }}_{=: c_*(f^1,f^2,\beta,T)} 
        \Norm{ \delta \zeta }^2_{ L^2 }
    +
    \underbrace{ c_2(\beta) \frac{e^{\beta K^* T }}{K_*} }_{=: c_{**}(f^1,f^2,\beta,T)} 
        \Esp{ \int_0^T \Abs{\delta f_t^\mu(Y',Z',U')}  \diff t }.
\end{align*}

\end{proof}

\section{Forward-Backward SDEs with Random Environment and Jumps}
\label{Sec:FBSDEs}

Now we proceed to investigate the system of coupled FBSDEs in admissible set-ups. We start by stating several preliminary results,  which are useful for proving the main theorem. Hereafter we assume that
\begin{align*}
    \Parentheses{ B, F , \Sigma, \Gamma, \zeta } 
    \in 
    \HH^2(W,\tilde{\eta}^\lambda) \times \HH^2(W,\tilde{\eta}^\lambda) \times \HH^2(W) 
        \times \HH^2(\tilde{\eta}^\lambda) \times L^2(\Ff_T; \RR^n)
\end{align*}
are arbitrary integrable processes, unless specified otherwise. 

\begin{remark}
\label{Lemma:Uncoupled-FBSDEs}
Note that, based on the previous results for uncoupled systems, the following FBSDEs have a unique solution $
    (X,Y,Z,U,M) 
    \in
    \Ss^2(\RR^d) \times \Ss^2(\RR^n) \times \HH^{2}(W) \times \HH^{2}(\tilde{\eta}^\lambda) \times 
        \Hh^{2,\perp}
$:

\begin{enumerate}[(i)]
    
    \item 
        Uncoupled (linear) forward component:
		\begin{align}
            \label{Eq:Uncoupled-Forward-Component-I}
            X_t
            =& 
            X_0 + \int_0^t B_s \diff s + \int_0^t \Sigma_s \diff W_s 
                + \int_{[0,t] \times \mathbf{R}} \Gamma_s(r) \tilde{\eta}^\lambda (\diff s, \diff r),
            \\ \nonumber
            Y_t
            =& 
            \Parentheses{ G X_T + \zeta }
            +
            \int_t^T \Parentheses{ G X_s + F_s } \diff s
            -
            \int_t^T Z_{s} \diff W_s
            -
            \int_{(t,T] \times \mathbf{R}} U_{s}(r) \tilde{\eta}^\lambda(\diff s, \diff r) 
            \\ \label{Eq:Uncoupled-Forward-Component-II}
            &-
            (M_T - M_t).
		\end{align}
		
    \item
        Uncoupled (linear) backward component:
		\begin{align}
            \label{Eq:Uncoupled-Backward-Component-I}
            X_t
            =& 
            X_0 + \int_0^t \Parentheses{ B_s + G^\top Y_s } \diff s 
                + \int_0^t \Parentheses{ \Sigma_s + G^\top Z_s } \diff W_s
                \\ \nonumber
                &+
                \int_{[0,t] \times \mathbf{R}} \Parentheses{ \Gamma_s(r) + G^\top U_s(r) } 
                    \tilde{\eta}^\lambda (\diff s, \diff r),
		  \\ \label{Eq:Uncoupled-Backward-Component-II}
            Y_t
            =& 
            \zeta + \int_t^T F_s \diff s - \int_t^T Z_{s} \diff W_s 
                - 
                \int_{(t,T] \times \RR_+} U_{s}(r) \tilde{\eta}^\lambda(\diff s, \diff r)
                -
                (M_T - M_t).
		\end{align}

\end{enumerate}

\end{remark}

\begin{lemma}
\label{Lemma:Expected-value-<GX,Y>}
Let $X \in \Ss^2(\RR^d)$ and $
    (Y,Z,U,M) \in \Ss^2(\RR^n) \times \HH^{2}(W) \times \HH^{2}(N^\lambda) \times \Hh^{2,\perp}
$ be the processes defined as
\begin{align*}
    X_t 
    &= 
    X_0 + \int_0^t B_s \diff s + \int_0^t \Sigma_s \diff W_s 
        + \int_{[0,t] \times \mathbf{R} } \Gamma_s(r) \tilde{\eta}^\lambda (\diff s, \diff r),
    \\
    Y_t
    &= 
    \zeta + \int_t^T F_s \diff s - \int_t^T Z_{s} \diff W_s 
        - \int_{(t,T] \times \mathbf{R} } U_{s}(r) \tilde{\eta}^\lambda(\diff s, \diff r) - (M_T - M_t).
\end{align*}
Then,
\begin{align}
	\label{Eq:Expected-value-<GX,Y>}
	\Esp{ ( X_T \cdot \xi )_G }
	&-
	\Esp{ ( X_0 \cdot Y_0 )_G }
	\\ \nonumber
	&=
	\Esp{ \int_0^T \Braces{
		- 
		( X_t \cdot F_t )_G
		+
		( B_t \cdot Y_t )_G
		+
		\Parentheses{ G\Sigma_t \cdot Z_t }_{ \mathrm{Fr} }
		+
		\Parentheses{ G \Gamma_t \cdot U_t }_{ \lambda_t }
	} \diff t }.
\end{align}

\end{lemma}

\begin{proof}
Equation \eqref{Eq:Expected-value-<GX,Y>} is obtained applying  Itô's formula and using  the orthogonality of $M$ with respect to $W$ and $\tilde{\eta}^\lambda$; namely,
\begin{align*}
    &\Esp{ ( X_T, Y_T )_G } - \Esp{ ( X_0, Y_0 )_G }
    =
    \sum_{i=1}^d \sum_{j=1}^n \Esp{ X^{(i)}_t G^{(i,j)} Y^{(j)}_t }
    \\
    &\qquad=
    \sum_{i=1}^d \sum_{j=1}^n G^{(i,j)} \Big\lbrace 
        \int_0^t \EE\Big[ 
            -X^{(i)}_s F^{(j)}_s + B^{(i)}_s Y^{(j)}_s
        \Big] \diff s
        +
        \Esp{ [X^{(i)}, Y^{(j)}]_t }
    \Big\rbrace.
\end{align*}
where the expected quadratic variation is given by
\begin{align*}
    \Esp{ [X^{(i)}, Y^{(j)}]_t }
    =
    \int_0^t \EE\Big[
        \Parentheses{ \Sigma^{(i,\cdot)}_s \cdot Z^{(j,\cdot)}_s }_\mathrm{Fr} 
        +
        \Parentheses{ \Gamma^{(i,\cdot)}_s \cdot U^{(j,\cdot)}_s }_{ \lambda_s }
    \Big] \diff s.
\end{align*}

\end{proof}

We are now ready to prove the uniqueness statement of our main theorem:
\begin{proof}[Proof of Theorem \ref{Thm:Existence-Uniqueness-FBSDE} (Uniqueness)]
Let $(X,Y,Z,U,M)$ and $(X',Y',Z',U',M')$ be two solutions of \eqref{Eq:Forward-Component}-\eqref{Eq:Initial-Terminal-Conditions}. Then, the difference is given by
\begin{align}
    \label{Eq:(FBSDE)delta-X}
    \delta{X}_t
    =&
    \int_0^t \delta{b}^\mu_s \diff s
    +
    \int_0^t \delta{\sigma}_s^\mu \diff W_s
    +
    \int_{[0,t] \times \mathbf{R} } \delta{\gamma}_s^\mu(r) \tilde{\eta}^\lambda ( \diff s, \diff r),
    \\ \label{Eq:(FBSDE)delta-Y}
    \delta{Y}_t 
    =&
    \delta{g}^\mu
    +
    \int_t^T \delta{f}^\mu_s \diff s
    -
    \int_t^T \delta{Z}_{s} \diff W_s
    -
    \int_{(t,T] \times \mathbf{R} } \delta U_s(r) \tilde{\eta}^\lambda ( \diff s, \diff r)
    \\ \nonumber
    &-
    \Parentheses{ \delta{M}_T - \delta{M}_t },
\end{align}
with the initial condition $\delta X_0 \equiv 0$.
Lemma \ref{Lemma:Expected-value-<GX,Y>} and the definition of $A^\mu$ in (\ref{DefA}) yield
\begin{align*}
    \Esp{ ( G \delta X_T \cdot \delta g^\mu ) }
    =
    \Esp{ \int_0^T \Parentheses{ \delta{A}_t^\mu \cdot \delta (X_t,Y_t,Z_t,U_t) } \diff t }.
\end{align*}
Then, from the $G$-mononicity Condition {\bf A}.6:
\begin{align}
	\nonumber
	0 
	\leq&
	- 
	\beta_1 \Esp{ \int_0^T \Abs{ G\delta{X}_{s} }^2 \diff s }
	- 
	\beta_2 \Esp{ \int_0^T \Braces{ 
		\Abs{ G^\top \delta{Y}_{s} }^2 
		+ 
		\Norm{ G^\top \delta{Z}_{s} }^2_{\mathrm{Fr}} 
		+ 
		\Norm{ G^\top \delta{U}_{s} }^2_{\lambda_t}
	} \diff s }
	\\ \label{Eq:(FBSDE)Majorant}
	&-
	\beta_3 \Esp{ \Abs{  G \delta{X}_T }^2 }.
\end{align}

When $d<n$, we have that $\beta_1,\beta_3 > 0$. From the previous inequality, $\Abs{ G \delta{X}_t }^2 = 0$ $\diff \PP \otimes \diff t$-a.s., which implies that $X$ and $X'$ are $\PP$-indistinguishable. In consequence, $Y_T = g^\mu(X_T)=g^\mu(X_T') = Y_T'$. Thus, from the uniqueness of solutions for BSDEs, it follows that $(Y,Z,U,M) = (Y',Z',U',M')$  $\PP$-a.s.

On the other hand, if $d\geq n$, we must have that $\beta_2 > 0$. From an analogous argument to the one given in \eqref{Eq:(FBSDE)Majorant} we get that $(Y,Z,U)=(Y',Z',U')$ $\diff \PP \otimes \diff t$-a.e.. Thus, by the uniqueness result of (forward) SDEs, $X = X'$ up to $\PP$-indistinguishably. In order to verify the equality of $M = M'$, from \eqref{Eq:(FBSDE)delta-Y} note that
\begin{align*}
    &\Norm{ \delta M }^2_{ \Hh^{2,\perp} }
    =
    \Esp{ \delta M^2_T }
    \\
    &=
    \Esp{ \Parentheses{
		\delta{g}^\mu
		+
		\int_0^T \delta{f}^\mu_s \diff s
		-
		\int_0^T \delta{Z}_{s} \diff W_s
		-
		\int_{ (0,T] \times \mathbf{R}} \delta{U}_{s}(r) \tilde{\eta}^\lambda( \diff s, \diff r)
		-
		\delta Y_0
    }^2 }
    =0.
\end{align*}

\end{proof}

In order to  complete the proof of  Theorem \ref{Thm:Existence-Uniqueness-FBSDE} -with the existence part-  the following results   are needed,  providing  the ground basis of the continuation argument. To ease our notation, the parameterized family   indexed by $\alpha\in[0,1]$ is written as
\begin{align*}
    (b,\sigma,\gamma,f)_s^\mu( X^\alpha, Y^\alpha, Z^\alpha, U^\alpha ) 
    = 
    (b,\sigma,\gamma,f)_s^{\mu,\alpha}.
\end{align*}

Notice that, depending on the relation between $n$ and $d$, each case is treated separately; see Remark \ref{dimen}. On the other hand, the  role of the processes $ B, F , \Sigma, \Gamma$  will become apparent in the iterative procedure proposed after the next results.
\begin{lemma}
\label{Lemma:Continuation-d}
Consider the following two sets of FBSDEs parametrized by $\alpha \in [0,1]$  and $(B, F , \Sigma, \Gamma, \zeta)$:
\begin{enumerate}[(i)]
    
    \item
        When $d < n$ (i.e. $\beta_1, \beta_3 > 0$),  let:
        \begin{align}
            \label{Eq:(FBSDE)Continuation-d-Forward}
            X_t^\alpha
            =& 
            X_0
            +
            \int_0^t \Braces{ \alpha b_s^{\mu,\alpha} + B_s } \diff s
            +
            \int_0^t \Braces{ \alpha \sigma_s^{\mu,\alpha} + \Sigma_s } \diff W_s
            \\ \nonumber
            &+
            \int_{[0,T] \times \mathbf{R}} \Braces{ \alpha \gamma_s^{\mu,\alpha}(r) + \Gamma_s(r) } 
                \tilde{\eta}^\lambda( \diff s, \diff r ),
            \\ \label{Eq:(FBSDE)Continuation-d-Backward}
            Y_t^\alpha
            =& 
            Y_T^\alpha
            +
            \int_t^T \Braces{ \alpha f_s^{\mu,\alpha} + F_s + (1 - \alpha) \beta_1 G X_s^\alpha } \diff s
            -
            \int_t^T Z_s^\alpha \diff W_s
            \\ \nonumber
            &-
            \int_{ (t,T] \times \mathbf{R} } U_s^\alpha(r) \tilde{\eta}^\lambda( \diff s, \diff r )
            -
            ( M_T^\alpha - M_t^\alpha ),
            \\ \label{Continuation-d-Initial-Terminal-Conditions}
            X_0 \in& L^2(\Ff_0;\RR^d),
            \quad 
            Y_T = \alpha g^\mu(X^\alpha_T) + \zeta + (1-\alpha)GX_T^\alpha,
        \end{align}
        
    \item
        When $d \geq n$ (i.e. $\beta_2 > 0$), define:  
        \begin{align}
            \label{Eq:(FBSDE)Continuation-n-Forward}
            X_t^\alpha
            =& 
            \int_0^t \Braces{ \alpha b_s^{\mu,\alpha} + B_s - (1 - \alpha) \beta_2 \Parentheses{ G^\top Y_s^\alpha} } 
                \diff s
            \\ \nonumber
            &+
            \int_0^t \Braces{ 
                \alpha \sigma_s^{\mu,\alpha} + \Sigma_s - (1 - \alpha) \beta_2 \Parentheses{ G^\top Z_s^\alpha} 
            } \diff W_s
            \\ \nonumber
            &+
            \int_{[0,t] \times \mathbf{R} } \Braces{ 
                \alpha \gamma_s^{\mu,\alpha} + \Gamma_s(r) - (1 - \alpha) \beta_2 \Parentheses{ G^\top U^\alpha_s(r) } 
            } \tilde{\eta}^\lambda ( \diff s, \diff r ),
            \\ \label{Eq:(FBSDE)Continuation-n-Backward}
            Y_t^\alpha
            =& 
            Y_T
            +
            \int_t^T \Braces{ \alpha f_s^{\mu,\alpha} + F_s } \diff s
            -
            \int_t^T Z_s^\alpha \diff W_s
            -
            \int_{(t,T] \times \mathbf{R}} U_s^\alpha(r) \tilde{\eta}^\lambda( \diff s, \diff r )
            \\ \nonumber
            &-
            ( M_T^\alpha - M_t^\alpha )
            \\ \label{Eq:(FBSDE)Continuation-n-Initial-Terminal-Conditions}
            X_0 \in& L^2(\Ff_0;\RR^d),
            \quad 
            Y_T = \alpha g^\mu(X^\alpha_T) + \zeta.
        \end{align}
\end{enumerate}

Assume that there exists $\alpha_0 \in [0,1)$ such that, for each $(B,F,\Sigma,\Gamma,\zeta)$, there exists a solution $( X^{\alpha_0}, Y^{\alpha_0}, Z^{\alpha_0}, U^{\alpha_0}, M^{\alpha_0} )$ of the corresponding FBSDEs in (i) or (ii). 
Then, there exists a positive constant $\epsilon_0 > 0$ such that, for each $\epsilon \in [0, \epsilon_0]$, with $\epsilon_0+\alpha_0\leq 1$, there exists a process $( X^{\alpha_0 + \epsilon}, Y^{\alpha_0 + \epsilon}, Z^{\alpha_0 + \epsilon}, U^{\alpha_0 + \epsilon}, M^{\alpha_0 + \epsilon} )$ solving the associated system  of FBSDEs 
with parameter $\alpha_0 + \epsilon$.
\end{lemma}

\begin{remark}
Note that by taking $\alpha = 1$ the existence of solutions to \eqref{Eq:(FBSDE)Continuation-d-Forward}-\eqref{Continuation-d-Initial-Terminal-Conditions} and \eqref{Eq:(FBSDE)Continuation-n-Forward}-\eqref{Eq:(FBSDE)Continuation-n-Initial-Terminal-Conditions} implies the existence of solutions to \eqref{Eq:Forward-Component}-\eqref{Eq:Initial-Terminal-Conditions}, as we can take $(B,F,\Sigma,\Gamma,\zeta)$ identically zero. Besides, observe that when $\alpha = 0$ the FBSDE reduces to an uncoupled system of the form \eqref{Eq:Uncoupled-Forward-Component-I}-\eqref{Eq:Uncoupled-Backward-Component-I} (resp. \eqref{Eq:Uncoupled-Forward-Component-II}-\eqref{Eq:Uncoupled-Backward-Component-II}), which  has a solution due to Lemma \ref{Lemma:Uncoupled-FBSDEs}. The next proof  will make clear that $\epsilon_0$ does not depend on $\alpha_0$, which becomes an important observation when the procedure is iterated.
\end{remark}
\begin{proof}
Let $\alpha_0$ be as in the hypothesis and  define the mapping $\Phi^{\alpha_0}$ as
\begin{align*}
    (B,F,\Sigma,\Gamma,\zeta) 
    \mapsto 
    \Phi^{\alpha_0}(B,F,\Sigma,\Gamma,\zeta) 
    := 
    ( X^{\alpha_0}, Y^{\alpha_0}, Z^{\alpha_0}, U^{\alpha_0}, M^{\alpha_0} ),
\end{align*}
where $( X^{\alpha_0}, Y^{\alpha_0}, Z^{\alpha_0}, U^{\alpha_0}, M^{\alpha_0} )$ is the solution process to \eqref{Eq:(FBSDE)Continuation-d-Forward}-\eqref{Continuation-d-Initial-Terminal-Conditions} associated with the inputs. Now, let
\begin{align}
    \label{Eq:(FBSDE)Preimage-of-Contraction-d}
    (x,y,z,u,m,x_T) 
    \in 
    \HH^2(W,\tilde{\eta}^\lambda) \times \HH^2(W,\tilde{\eta}^\lambda) \times \HH^{2}(W) 
        \times \HH^{2}(\tilde{\eta}^\lambda) \times \Hh^{2,\perp} \times L^2(\Ff_T).
\end{align}
Using the fact that $(B,F,\Sigma,\Gamma,\zeta)$ has been taken arbitrarily, Remark \ref{Lemma:Uncoupled-FBSDEs} implies that \eqref{Eq:(FBSDE)Continuation-d-Forward}-\eqref{Continuation-d-Initial-Terminal-Conditions}  has a unique solution process for $(B^\epsilon, F^\epsilon, \Sigma^\epsilon, \Gamma^\epsilon,\zeta^\epsilon)$, defined as
\begin{align*}
    &B_t^\epsilon  :=  B_t + \epsilon b_t^\mu(x,y,z,u),
    &
    &F_t^\epsilon := F_t + \epsilon \Parentheses{ f_t^\mu(x,y,z,u) - \beta_1 G x_t },
    \\
    &\Sigma_t^\epsilon := \Sigma_t + \epsilon \sigma_t^\mu(x,y,z,u),
    &
    &\zeta^\epsilon := \zeta + \epsilon \Parentheses{ g^\mu(x_T) - G x_T },
    \\
    &\Gamma_t^\epsilon := \Gamma_t + \epsilon \gamma_t^\mu(x,y,z,u),
\end{align*}
with $\epsilon > 0$. Composing the previous maps, let $\Phi^{\alpha_0, \epsilon}$ be defined by 
\begin{align}
    \label{Eq:(FBSDE)Mapping-d}
    (x,y,z,u,m,x_T) 
    \mapsto 
    \Phi^{\alpha_0, \epsilon}(x,y,z,u,x_T) 
    :&=
    \Phi^{\alpha_0}(B^\epsilon, F^\epsilon, \Sigma^\epsilon, \Gamma^\epsilon,\zeta^\epsilon).
\end{align}
Notice that a fix point of $\Phi^{\alpha_0, \epsilon}$ is a solution of  the 
FBSDE system with parameter $\alpha_0+\epsilon$.
Hence, we proceed to prove that, for each $\epsilon \in [0,\epsilon_0]$, this mapping is contractive, with $\epsilon_0 > 0$ to be specified later.  

First, let us precise that the domain and range of this mapping is given by
\begin{align*}
    \Hh^2(\RR^d) \times \Hh^2(\RR^n) \times \HH^{2}(W) \times \HH^{2}(\tilde{\eta}^\lambda) \times \Parentheses{ \Ll^2(W) \oplus \Ll^2(\tilde{\eta}^\lambda)} ^{\perp}.
\end{align*}
Note that this is a Banach space with respect to a norm equivalent to
\begin{align*}
    \Norm{\cdot}_* 
    = \Norm{\cdot}_{\Ss^2(\RR^d)} + \Norm{\cdot}_{\Ss^2(\RR^n)} + \Norm{\cdot}_{\HH^2(W)} 
        + \Norm{\cdot}_{\HH^2(\tilde{\eta}^\lambda)} + \Norm{\cdot}_{\Hh^{2}};
\end{align*}
see Theorems IV.1 and V.2 from \cite{protter_stochastic_2005},  

We now write
\begin{align*}
    (X,Y,Z,U,M,X_T) := \Phi^{\alpha_0,\epsilon} (x,y,z,u,m,x_T);
\end{align*}
similarly, for another process $(x',y',z',u',m',x_T')$ in the same space as \eqref{Eq:(FBSDE)Preimage-of-Contraction-d}, we write 
\begin{align*}
    (X',Y',Z',U',M',X_T') := \Phi^{\alpha_0,\epsilon} (x',y',z',u',m',x_T').
\end{align*}

Note that the difference process satisfies
\begin{align}
    \label{Eq:(FBSDE)delta-X-d}
    \delta{X}_t 
    =&
    \int_0^t \Braces{ \alpha_0 \delta{b}^\mu_s(X,Y,Z,U) + \epsilon \delta{b}^\mu_s(x,y,z,u)} \diff s
    \\ \nonumber
    &+
    \int_0^t \Braces{ \alpha_0 \delta{\sigma}^\mu_s(X,Y,Z,U) + \epsilon \delta{\sigma}^\mu_s(x,y,z,u) }\diff W_s
    \\ \nonumber
    &+
    \int_{[0,t] \times \mathbf{R}} \Braces{ 
        \alpha_0 \delta{\gamma}^\mu_s(V,U) + \epsilon \delta{\gamma}^\mu_s(v,u) 
    } \tilde{\eta}^\lambda( \diff s, \diff r ),
    \\ \label{Eq:(FBSDE)delta-Y-d}
    \delta{Y}_t 
    =&
    \Braces{ 
        \alpha_0 \delta{g}^\mu( X_T ) + (1 - \alpha_0) G \delta{X}_T 
        + 
        \epsilon \Brackets{ \delta{g}^\mu( {x}_T ) - G \delta{x}_T } 
    }
    \\ \nonumber
    &+
    \int_t^T \Braces{ 
        \alpha_0 \delta{f}^\mu_s(X,Y,Z,U) 
        + 
        (1-\alpha_0) \beta_1 G \delta{X}_s 
        +
		\epsilon \Brackets{ \delta{f}^\mu_s(x,y,z,u) - \beta_1 G \delta{x}_s}
    } \diff s
    \\ \nonumber
    &-
    \int_t^T \delta{Z}_{s} \diff W_s
    -
    \int_{(t,T] \times \mathbf{R} } \delta{U}_{s}(r) \tilde{\eta}^\lambda( \diff s, \diff r )
    -
    (\delta M_T - \delta M_t).
\end{align}

Let $q > 0$. From Lemma \ref{Lemma:Expected-value-<GX,Y>} we obtain
\begin{align*}
    e^{qt} \Esp{ ( \delta X_t \cdot \delta Y_t )_G }
    =&
    q \int_0^t e^{rs} \Esp{ ( \delta X_s \cdot \delta Y_s )_G } \diff s
    -
    (1-\alpha_0) \beta_1 \int_0^t e^{qs} \Esp{ \Abs{ G \delta{X}_s }^2 } \diff s
    \\
    &+
    \epsilon \beta_1 \int_0^t e^{qs} \Esp{ \Parentheses{ \delta X_s \cdot G \delta{x}_s }_G } \diff s
    \\
    &+
    \alpha_0 \int_0^t e^{qs} \Esp{ \Parentheses{ \delta{A}^\mu_s(X,Y,Z,U) \cdot \delta (X,Y,Z,U)_s } } \diff s
    \\
    &+
    \epsilon \int_0^t e^{qs} \Esp{ \Parentheses{ \delta A^\mu_s(x,y,z,u) \cdot \delta (X,Y,Z,U)_s } } \diff s.
\end{align*}
Moreover, from the $G$-monotonicity and Lipschitz conditions,
\begin{align*}
    e^{qt} &\Esp{ ( \delta X_t \cdot \delta Y_t )_G }
    \\
    \leq&
    q \int_0^t e^{qs} \Esp{ ( \delta X_s \cdot \delta Y_s )_G } \diff s
    -
    (1-\alpha_0) \beta_1 \int_0^t e^{qs} \Esp{ \Abs{ G \delta{X}_s }^2 } \diff s
    \\
    &+
    \frac{\epsilon}{2} \beta_1 \int_0^t e^{qs} \Esp{ \Abs{ G \delta X_s }^2 + \Abs{ G \delta{x}_s }^2 } \diff s
    \\ 
    &+
    \alpha_0 \Esp{ \int_0^t e^{qs} \Braces{
		-\beta_1 \Abs{ G \delta X_s }^2
		- 
		\beta_2 \Parentheses{ 
            \Abs{ G^\top \delta{Y}_{s} }^2 
            + 
            \Norm{ G^\top \delta{Z}_{s} }^2_{\mathrm{Fr}}
            + 
            \Norm{ G^\top \delta{U}_{s-} }^2_{\lambda_t}
        }
    } \diff s }
    \\ 
    &+
    \frac{\epsilon}{2} c_K \int_0^t e^{qs} \Esp{
        \Abs{ \delta x_s }^2 
        +
        \Abs{ \delta y_s }^2 
        +
        \Norm{ \delta z_s }^2_{ \mathrm{Fr} }
        + 
        \Norm{ \delta u_s }^2_{\lambda_s}
    } \diff s
    \\ 
    &+
    \frac{\epsilon}{2} \int_0^t e^{qs} \Esp{
        \Abs{ G \delta X_s }^2 
		+
		\Abs{ G^\top \delta Y_s }^2 
		+
		\Norm{ G^\top \delta Z_s }^2_{ \mathrm{Fr} }
		+ 
		\Norm{ G^\top \delta U_s }^2_{\lambda_s}
    } \diff s
    \\
    \leq&
    q \int_0^t e^{qs} \Esp{ ( \delta X_s \cdot \delta Y_s )_G } \diff s
    -
    \beta_1 \int_0^t e^{qs} \Esp{ \Abs{ G \delta{X}_s }^2 } \diff s
    \\ 
    &+
    \alpha_0 \Esp{ \int_0^t e^{qs} \Braces{
		- \beta_2 \Parentheses{ 
            \Abs{ G^\top \delta{Y}_{s} }^2 
            + 
            \Norm{ G^\top \delta{Z}_{s} }^2_{\mathrm{Fr}}
            + 
            \Norm{ G^\top \delta{U}_{s} }^2_{\lambda_t}
		}
    } \diff s }
    \\ 
    &+
    \frac{\epsilon}{2} (c_K \vee \beta_1 c_G ) \int_0^t e^{qs} \Esp{
        \Abs{ \delta x_s }^2 
        +
        \Abs{ \delta y_s }^2
        +
        \Norm{ \delta z_s }^2_{ \mathrm{Fr} }
        + 
        \Norm{ \delta u_s }^2_{\lambda_s}
    } \diff s 
    \\ 
    &+
    \frac{\epsilon}{2} (1 \vee \beta_1) c_G \int_0^t e^{qs} \Esp{
        \Abs{ \delta X_s }^2 
		+
		\Abs{ \delta Y_s }^2
		+
		\Norm{ \delta Z_s }^2_{ \mathrm{Fr} }
		+ 
		\Norm{ \delta U_s }^2_{\lambda_s}
    } \diff s.
\end{align*}
Then,  Gr\"{o}nwall's inequality yields,
\begin{align}
    \nonumber
    &e^{qt} \Esp{ ( \delta X_t \cdot \delta Y_t )_G }
    +
    \beta_1 \int_0^t e^{qs} \Esp{ \Abs{ G \delta{X}_s }^2 } \diff s
    \\ \nonumber
    &\qquad+
    \alpha_0 \beta_2 \Esp{ \int_0^t e^{qs} \Braces{
        \Abs{ G^\top \delta{Y}_{s} }^2 
		+ 
		\Norm{ G^\top \delta{Z}_{s} }^2_{\mathrm{Fr}}
		+ 
		\Norm{ G^\top \delta{U}_{s} }^2_{\lambda_t}
    } \diff s }
    \\ \label{Eq:(FBSDE)Gronwall-d}
    &\leq
    \epsilon \frac{c_K \vee c_G \vee \beta_1 c_G }{2} e^{qt} \left\lbrace
		\int_0^t e^{qs} \Esp{
            \Abs{ \delta x_s }^2
            + 
            \Abs{ \delta y_s }^2
            +
            \Norm{ \delta z_s }^2_{ \mathrm{Fr} }
            + 
            \Norm{ \delta u_s }^2_{\lambda_s}
        } \diff s 
        \right.
        \\ \nonumber
        &\qquad\qquad\qquad\qquad\qquad\qquad
        \left.
		+
		\int_0^t e^{qs} \Esp{
		  \Abs{ \delta X_s }^2
            + 
            \Abs{ \delta Y_s }^2
            +
            \Norm{ \delta Z_s }^2_{ \mathrm{Fr} }
            + 
            \Norm{ \delta U_s }^2_{\lambda_s}
		} \diff s 
	\right\rbrace,
\end{align}
for all $t \in [0,T]$. On the other hand, due to the terminal conditions,
\begin{align}
    \nonumber
    \Esp{ e^{qT} \Parentheses{ \delta{X}_T \cdot \delta Y_T }_G }
    =&
    \alpha_0 \Esp{ e^{qT} \Parentheses{ 
		\delta X_T 
		\cdot 
		\delta{g}^\mu( X_T ) 
    }_G }
    +
    (1- \alpha_0) \Esp{ e^{qT} \Abs{ G \delta{X}_T }^2 }
    \\ \nonumber
    &+
    \epsilon \Esp{ e^{qT} \Braces{
		\Parentheses{ \delta X_T, \delta{g}^\mu( x_T ) }_G
		-
		\Parentheses{ G \delta{X}_T, G \delta{x}_T }
    } }
    \\ \label{Eq:(FBSDE)Terminal-d}
    \geq&
    \Brackets{ \alpha_0 \beta_3 + (1- \alpha_0) } \Esp{ e^{qT} \Abs{ G \delta{X}_T }^2 }
    \\ \nonumber
    &-
    \epsilon \Esp{ e^{qT} \Braces{
		\frac{ \Abs{ G \delta{X}_T }^2 + \Abs{ G \delta{x}_T }^2 }{2}
		-
		\Parentheses{ \delta X_T \cdot \delta{g}^\mu( {x}_T ) }_G
  } }.
\end{align}
Note that $\alpha_0 \beta_3 + (1- \alpha_0) > 1 \wedge \beta_3$. Combining \eqref{Eq:(FBSDE)Gronwall-d} and \eqref{Eq:(FBSDE)Terminal-d} and letting $q \to 0$ yields
\begin{align}
	\nonumber
	&\underbrace{(1 \wedge \beta_1 \wedge \beta_3 )}_{ > 0 }\Braces{
		\Norm{ \delta X }^2 
		+
		\Norm{ \delta{X}_T }^2_{ L^2 }
	}
	\\ \label{Eq:(FBSDE)Forward-Bound-d}
	&\quad\leq 
	\beta_1 \Norm{ \delta X }^2 
	+
	\alpha_0 \beta_2 \Norm{ \delta (Y,Z,U) }^2
	+
	\Brackets{ \alpha_0 \beta_3 + (1- \alpha_0) } \Norm{ \delta{X}_T }^2_{ L^2 }
	\\ \nonumber
	&\quad\leq
	\epsilon \frac{ c_K \vee c_G \vee \beta_1 c_G }{2 c_G } \Braces{
		\Norm{ \delta(x,y,z,u) }^2
		+
		\Norm{ \delta x_T }_{ L^2 }^2
		+
		\Norm{ \delta(X,Y,Z,U) }^2
		+
		\Norm{ \delta X_T }_{ L^2 }^2
	}.
\end{align}
Now, we make use of the \emph{a priori} estimates obtained in Lemma \ref{Lemma:(BSDE)A-priori-estimates}, as well as the fact that $\alpha_0 < 1$:
\begin{align}
	\label{Eq:(FBSDE)a-priori-d}
	\Norm{ \delta (Y,Z,U,M) }^2
	\leq
	c_f \Norm{ \delta Y_T }^2
	\leq
	\underbrace{ 
		c_f \Parentheses{ \underline{\underline{K}} + c_G } 
	}_{ 
		=: c_{ f, g, G } 
	} \Parentheses{ 
		\Norm{ \delta X_T }^2_{ L^2 } 
		+ 
		\epsilon^2 \Norm{ \delta x_T }^2_{ L^2 } 
	}.
\end{align}
Consequently, by setting $
    \underline{c} 
    := 
    (2 c_G)^{-1}( 1 \wedge \beta_1 \wedge \beta_3 )^{-1} (c_{f,g,G} + 1) (c_K \vee c_G \vee \beta_1 c_G)
$,
\begin{align*}
    &\Norm{ \delta \Phi^{\alpha_0,\epsilon}(x,y,z,u,m,x_T)}^2
    =
    \Norm{ \delta (X, Y, Z, U, M, {X}_T) }^2
    \\
    &\quad=
    \Norm{ \delta (Y,Z,U,M) }^2 + \Norm{ \delta{X} }^2
    +
    \Norm{ \delta{X}_T }_{ L^2 }^2 
    \\ (\mbox{eq. } \eqref{Eq:(FBSDE)a-priori-d})
    &\quad\leq
    ( c_{f,g,G,\alpha_0} + 1 ) \Braces{ \Norm{ \delta X_T }^2_{ L^2 } + \Norm{ \delta X }^2 }
    +
    c_{f,g,G} \epsilon^2 \Norm{ \delta x_T }^2_{ L^2 }
    \\ (\mbox{eq. } \eqref{Eq:(FBSDE)Forward-Bound-d})
    &\quad\leq
    \underline{c} \epsilon \Braces{ \Norm{ \delta (x,y,z,u,m,x_T) }^2 + \Norm{ \delta (X,Y,Z,U,M,X_T) }^2 }
    +
    c_{f,g,G} \epsilon^2 \Norm{ \delta x_T }^2_{ L^2 }.
\end{align*}

Observe that by taking $\epsilon = \Parentheses{ 2 n_0 \underline{c} }^{-1}$, with $n_0 \in \NN$, we get $c_{f,g,G} \epsilon^2 \leq n_0^{-2}$. Hence, for $n_0 > 1$,
\begin{align*}
	\Norm{ \delta (X,Y,Z,U,M,X_T) }^2
	\leq
	\frac{2}{n_0-1}
	\Norm{ \delta (x,y,z,u,m,x_T) }^2.
\end{align*}
Thus, by taking $\epsilon_0 := \Parentheses{ 8 \underline{c} }^{-1}$, the mapping $\Phi^{\alpha_0, \epsilon}$ in \eqref{Eq:(FBSDE)Mapping-d} is a contraction  for  $\epsilon \in [0,\epsilon_0]$.

The proof for the case $d \geq n$ follows the same procedure, with a slight modification in the definition of the  processes $(B^\epsilon, F^\epsilon, \Sigma^\epsilon, \Gamma^\epsilon,\zeta^\epsilon)$, taking instead
\begin{align*}
    &B_t^\epsilon  :=  B_t + \epsilon \Parentheses{ b_t^\mu(x,y,z,u) + \beta_2 G^\top y_t },
    &
    &F_t^\epsilon := F_t + \epsilon f_t^\mu(x,y,z,u),
    \\
    &\Sigma_t^\epsilon := \Sigma_t + \epsilon \Parentheses{ \sigma_t^\mu(x,y,z,u) + \beta_2 G^\top z_t },
    &
    &\zeta^\epsilon := \zeta + \epsilon g^\mu(x_T),
    \\
    &\Gamma_t^\epsilon := \Gamma_t + \epsilon \Parentheses{ \gamma_t^\mu(x,y,z,u) + \beta_2 G^\top u_t }.
\end{align*} 
\end{proof}

\begin{proof}[Proof of Theorem \ref{Thm:Existence-Uniqueness-FBSDE} (Existence)]

We present the proof only for the case  $d < n$, since the arguments to prove  the other one ($d \geq n$)  are similar and can be adapted by making the corresponding changes.  From Remark \ref{Lemma:Uncoupled-FBSDEs} it follows that there exists a solution to \eqref{Eq:(FBSDE)Continuation-d-Forward}-\eqref{Continuation-d-Initial-Terminal-Conditions}
when $\alpha_0 = 0$ and, by Lemma \ref{Lemma:Continuation-d} part (i),  there exists a constant $\epsilon_0 > 0$   such that, for each $\epsilon \in [0, \epsilon_0]$, there exists a unique solution $(X^{ \epsilon}, Y^{ \epsilon}, Z^{ \epsilon}, U^{ \epsilon}, M^{ \epsilon})$ of \eqref{Eq:(FBSDE)Continuation-d-Forward}-\eqref{Continuation-d-Initial-Terminal-Conditions}. Now  iterate this procedure  $k$ steps, with $\alpha_{i+1}:=\alpha_i + \epsilon_0$, and select $k$ in such a way that $k \epsilon_0 \leq  1<  (1 +k) \epsilon_0$, concluding  that \eqref{Eq:(FBSDE)Continuation-d-Forward}-\eqref{Continuation-d-Initial-Terminal-Conditions}   has a unique solution with $\alpha = 1$. Finally, since the initial selection of $\Parentheses{ B, F , \Sigma, \Gamma, \zeta }$ was arbitrary, we get the existence of solution to (\ref{Eq:Forward-Component})-(\ref{Eq:Initial-Terminal-Conditions}).

\end{proof}

\section{Some Examples}
\label{Sec:SE}
In this section we present some examples to illustrate the presence of an exogenous random environment, describing its role in the dynamics of the system.   We will begin by examining the types of common noise that fit the proposed definitions.

\begin{example}
Given a random environment $\mu$ and an $\FF^0\otimes\Rr$-adapted simple point process $\eta$, let $\lambda : \Omega^0 \times [0,T] \times \Rr \times \Pp_2(\RR^d) \times \Mm_c^*([0,T] \times \mathbf{R}) \to \RR_+ $ be a measurable function of the form
\begin{align}
    \label{Eq:Additive-Kernel}
    \lambda( \omega^0, t, C, \nu, \eta )
    =
    \psi^0( \omega^0, t, C )
    +
    \psi^1( \nu, C )
    +
    \int_{ [0,t) \times (C \cap \mathbf{R}) } \psi^2(\omega^0, s, r') \eta(\omega^0, \diff s, \diff r'),
\end{align}
where 
\begin{enumerate}
    \item 
        $\psi^0 : \Omega \times [0,T] \times \Rr \to \RR_+$ is a locally integrable kernel from $(\Omega^0 \times [0,T], \Ff^0 \otimes \Bb([0,T])$ into $(\mathbf{R},\Rr)$; see (\ref{Eq:(Omega,Ff,PP)}), independent of $\mu$ and $\eta$.
        
    \item 
        $\psi^1 : \Pp_2(\RR^d) \times \Rr \to \RR_+$ is a measurable map such that, for each $C \in \Rr$, $\psi^1(\cdot,C)$ is a  non-negative Lipschitz function with respect to the metric $W_2$   and, for each $\nu \in \Pp_2(\RR^d)$, $\psi^1(\nu,\cdot)$ is a $\sigma$-finite measure on $(\mathbf{R},\Rr)$.
       
    \item 
        $\psi^2 \in \HH^2(\tilde{\eta})$.
\end{enumerate}
Then, $K^\lambda(\omega^0,t,C) := \lambda( \omega^0, t, C, \mu_{t-}, \eta^\lambda )$ defines an admissible $\FF^{0} \vee \FF^{\psi^0}$-intensity kernel candidate.
\end{example}

Note that the form of the intensity function $\lambda$ in \eqref{Eq:Additive-Kernel} allows us to consider a wide range of point processes as candidates for admissible noise, including Cox, semi-Cox \cite{bremaud_point_2020} and Hawkes processes, as long as we are able to work with flows on the Wasserstein space of probability measures. This is done with relative ease for mappings $\psi^1(\nu,\cdot) = \psi^1(\nu)$ admitting the integral representation as 
\begin{align}
    \label{Eq:tilde-psi}
    \psi^1(\nu) = \int_\Xx \tilde{\psi}^1(x) \nu(\diff x),
\end{align}
for some Lipschitz function $\tilde{\psi}^1$. In this case, an application of Jensen's inequality and the definition of the $W_2$-distance show that $\psi^1$ is Lipschitz on $\Pp_2(\Xx)$. Hereafter we will refer to kernels  satisfying \eqref{Eq:Additive-Kernel} and \eqref{Eq:tilde-psi} as \textit{additive kernel} or simply as \textit{additive intensity} when the intensity function $\lambda$ in \eqref{Eq:Additive-Kernel} does not depend on $C$.

Borrowing a well-known example from \textit{Optimal Control} and \textit{Game Theory}, we present an application of our main results to the LQ model.
\begin{example}[One-Dimensional Linear System]
\label{Example:Toy-model-III}
For easiness of notation, we consider of a one-dimensional, unmarked common noise. 
Suppose our standing assumptions hold, with $N$ a Poisson process on $[0,T] \times \RR \times \RR_+$ with intensity measure $\diff t \otimes \delta_{ \{ 1 \} }(\diff r) \otimes \diff s$. By definition, there exists a sequence of random points in $[0,T] \times \RR \times \RR_+$, denoted by $\Braces{(\tau_n^*,1,\xi^*_n), n \geq 1}$ such that $N(\diff t, \diff r, \diff s) = \sum_n \delta_{(\tau_n^*,1,\xi^*_n)}(\diff t, \diff r, \diff s)$.

For each $n\geq 1$, let $\xi_n \equiv 1$  and $\tau_{0}^{\lambda*} := \tau_0^* \equiv 0$. Let
\begin{align*}
    &\lambda_0^*(\omega^0,t) := \lambda_0 + \psi^1 v_{r}( \mu_{t-}(\omega^0) ),
    &
    &\forall (\omega^0,t) \in \Omega^0 \times (0,T],
\end{align*}
where $\lambda_{0}, \psi^1$ are positive constants, and 
\begin{align*}
    &v_{r}( \nu ) := \int_{ B_r(0) } x^2 \nu( \diff x ),
    &
    &\forall \nu \in \Pp_2(\RR),
\end{align*}
where $B_r(0)$ is the ball of radius $r$ centered at the origin. Let $\psi^2:[0,T] \to \RR_+$ be a continuous function and define the intensity process recursively as follows. For each $(\omega^0,t) \in \Omega^0 \times [0,T]$,
\begin{align*}
    \lambda_{n+1}^*(\omega^0,t) - \lambda_{n}^*(\omega^0,t) 
    = 
    \begin{cases}
        \psi^2(t - \tau^{*}_{n+1}(\omega^0))\Ind{ (\tau^{*}_{n+1}(\omega^0),T] },   &\mbox{if } \lambda_n^*(\omega^0,\tau_{n+1}^*(\omega^0)) \geq \xi_{n+1}^*(\omega^0),
        \\
        0,  &\mbox{otherwise};
    \end{cases}
\end{align*}
for each $n \geq 1$. Then, define
\begin{align*}
    \lambda(\omega^0,t) :&= \lim_{n \to \infty}\lambda_{n}^*(\omega^0,t).
\end{align*}
Consequently, the admissible common noise is given by $(N^\lambda,1)$; i.e. the marked Poisson process with marks $\xi_n \equiv 1$ and jumping times $\Braces{ \tau_n^\lambda }$ given by
\begin{align*}
    \tau_{n+1}^{\lambda*} (\omega^0) 
    :&= 
    \begin{cases}
        \tau_{n+1}^*(\omega^0), &\mbox{if } \lambda_n^*(\omega^0,\tau_{n+1}^*(\omega^0)) \geq \xi_{n+1}^*(\omega^0),
        \\
        \tau_{n}^{\lambda*}(\omega^0), &\mbox{otherwise};
    \end{cases}
    \\
    \tau_{n+1}^{\lambda} (\omega^0) :&= \min_{ k \geq 0 } \Braces{ \tau_{k}^{\lambda*}(\omega^0) : \tau_{k}^{\lambda*}(\omega^0) > \tau_{n}^{\lambda} (\omega^0) }.
\end{align*}
Then, by Theorem \ref{Thm:Existence-of-admissible-settings}, there exists an admissible set-up $(\Omega,\Ff,\FF,\PP,X_0,\mu,\lambda,N^\lambda,W)$, such that $N^\lambda$ is a self-exciting process with $\FF$-stochastic intensity
\begin{align}
    \label{Eq:LQ-lambda}
    \lambda_t 
    = 
    \lambda_0 + \psi^1 v_r(\mu_{t-}) + \int_{[0,t)} \psi^2(t-s) \diff N^\lambda_s.
\end{align}

Now, let $b,f,\sigma,\gamma \in \RR \setminus \{ 0 \}$ and $f_1,f_2,g > 0$ be constants such that
\begin{align*}
    &\Abs{b} f_2 = 2
    &
    &\mbox{and}
    &
    &f_1 f_2 > \frac{f^2}{2} + 1.
\end{align*}
Let $\widehat{b}:= b - \frac{f}{f_2}$ and $\widehat{f}:= \frac{f_1}{2} - \frac{f^2}{f_2}$, and define the linear mapping $m:\Pp_2(\RR) \to \RR$ as
\begin{align*}
    &m( \nu ) := \int_{ \RR } x \nu( \diff x ),
    &
    &\forall \nu \in \Pp_2(\RR).
\end{align*}

Moreover, we define the coefficients of the forward component as
\begin{align}
    \label{Eq:LQ-b,sigma,gamma}
    &b_t^\mu(x,y) 
    = 
    \widehat{b} \Parentheses{ x - m(\mu_{t-}) } - \frac{1}{f_2}y,
    &
    &\sigma_t^\mu(x,y)  = \sigma,
    &
    &\gamma_t^\mu(x,y) = \gamma,
\end{align}
for all $x,y \in \RR$. On the other hand, we define the elements of the backward component as
\begin{align}
    \label{Eq:LQ-f,g}
    &f_t^\mu(x,y) = \widehat{f}(x - m(\mu_{t-}) ) + \widehat{b} y,
    &
    &g^\mu_T(x) = g \Parentheses{ x - m(\mu_{T-}) }.
\end{align}
In other words, we have coupled linear FBSDE
\begin{align}
    \label{Eq:LQ-X}
    X_t 
    =& 
    X_0 
    + 
    \int_0^t \Braces{ \widehat{b} \Parentheses{ X_{s} - m(\mu_{s-}) } - \frac{1}{f_2} Y_s } \diff s 
    + 
    \sigma W_t 
    + 
    \gamma \tilde{N}^\lambda_t
    \\ 
    \label{Eq:LQ-Y}
    Y_t 
    =& 
    g( X_T - m(\mu_{T-}) ) 
    + 
    \int_t^T \Braces{ \widehat{f} \Parentheses{ X_s - m( \mu_{s-} ) } + \widehat{b} Y_s } \diff s 
    - 
    \int_t^T Z_s \diff W_s 
    \\ \nonumber
    &- 
    \int_{t+}^T U_{s} \diff \tilde{N}^\lambda_s 
    - 
    ( M_T - M_t).
\end{align}

It can be verified easily that the coefficients in \eqref{Eq:LQ-b,sigma,gamma}-\eqref{Eq:LQ-f,g} satisfy the conditions from Theorem \ref{Thm:Existence-Uniqueness-FBSDE}, taking $G=1$ and $\beta_1 = f_1 - \frac{f^2 - 1}{f_2} > 0$, $\beta_2 = 0$, $\beta_3 = g>0$ as constants in the $G$-monotonicity condition. Therefore, there exists a unique process
\begin{align*}
    (X,Y,Z,U,M) 
    \in 
    \Hh^2(\RR) \times \Hh^2(\RR) \times \HH^{2}(W) \times \HH^{2}(\tilde{N}^\lambda) 
        \times \Hh^{2,\perp}(\RR)
\end{align*}
 satisfying \eqref{Eq:LQ-X}-\eqref{Eq:LQ-Y}, subject to the fixed random environment $\mu$ and Common Noise with stochastic intensity $\lambda$ given by \eqref{Eq:LQ-lambda}.
\end{example}

\begin{remark}[Quadratic Hamiltonian]
In a forthcoming paper, we will analyze FBSDEs motivaded by Example \ref{Example:Toy-model-III}, which are   related to  Hamiltonian systems appearing in the \emph{Linear-Quadratic Mean-Field Game with Common Poissonian Noise} characterised by the payoff
\begin{align}
    \label{Eq:LQ-Payoff}
    \EE\Big[
        &\int_0^T \Big\lbrace 
            \frac{f}{2} \Parentheses{ X_t - \Esp[^1]{X_{t-}} }^2
            +
            f_1 \Parentheses{ X_t - \Esp[^1]{X_{t-}} } \alpha_t
            +
            \frac{f_2}{2} \alpha_t^2
		\Big\rbrace \diff t 
        +
        \frac{g}{2} \Parentheses{ X_T - \Esp[^1]{X_{T-}} }^2
    \Big],
\end{align}
defined on the admissible set-up $(\Omega,\Ff,\FF,\PP,X_0,\Ll^1(X),\lambda,(N^\lambda,1),W)$, where $\Ll^1(\ \cdot\ ) = \Ll(\ \cdot\ \vert \FF^{X_0,N^\lambda} )$ and $\alpha$ is an admissible control (i.e. a real-valued $\FF$-adapted process), subject to the dynamic evolution 
\begin{align}
    \label{Eq:LQ-Dynamic}
    X_t 
    =& 
    X_0 
    + 
    \int_0^t \Braces{ b \Parentheses{ X_{s} - \Esp[^1]{ X_{s-} } } + \alpha_s } \diff s 
    + 
    \sigma W_t 
    + 
    \gamma \tilde{N}^\lambda_t.
\end{align}
Namely, \eqref{Eq:LQ-X}-\eqref{Eq:LQ-Y} corresponds to the adjoint system of \eqref{Eq:LQ-Payoff}-\eqref{Eq:LQ-Dynamic} with $\widehat{\alpha}$ being the control  maximizing the Hamiltonian, i.e.
\begin{align*}
    &\widehat{\alpha}_t
    :=
    - 
    \frac{ f }{f_2} \Parentheses{ X_{t} - \Esp[^1]{ X_{t-} } } 
    - 
    \frac{ 1 }{f_2} Y_{t},
    &
    &\diff \PP \otimes \diff t-\mbox{a.e.},
\end{align*}
and $Y$ being the adjoint, first order variation of \eqref{Eq:LQ-Payoff} near the equilibrium.

\end{remark}

Another important family of admissible noise that does not necessarily fall under the additive class is the \textit{Environment-dependent Markov Chain} presented next. The construction is made following  \cite{shao_propagation_2022}.
\begin{example}
  Let $\Ee=\Braces{ 1, \ldots, n }$ be a finite set of states and, for each $\nu \in \Pp_2(\RR^d)$, let $Q(\nu)=\Parentheses{ Q^{(i,j)}(\nu); i,j \in \Ee }$ be a conservative, irreducible $Q$-matrix such that each $Q^{(i,j)}$ is a  mapping  on $\Pp_2(\RR^d)$, which is Lipschitz with respect to $W_2$. Furthermore, assume that the mapping 
\begin{align*}
	\Ee \times \Pp_2(\RR^d)
	\ni
	(i,\nu)
	\mapsto
	Q^{(i)}(\nu)
	:=
	\sum_{j \in \Ee} Q^{(i,j)}(\nu)
\end{align*}
is uniformly bounded, i.e. there exists a positive constant $H_0$ such that 
\begin{align*}
	H_0 := \sup_{ \nu \in \Pp_2(\RR^d) } \sup_{i \in \Ee} Q^{(i)}(\nu) < \infty.
\end{align*}
The set  $\Ee \times \Ee$ is ordered according to the lexicographic ordering \cite{yin_hybrid_2010}. For all $i,j \in \Ee$ with $i \neq j$ and $\nu \in \Pp_2(\RR^d)$, let $\Gamma^{(i,j)}$ be the consecutive   left-closed, right-open intervals of $\RR_+$, each one having length $Q^{(i,j)}(\nu)$, and $\Gamma^{(i,i)} = \emptyset$:
\begin{align*}
	\Gamma^{(i,j)}(\nu)
	:&=
	\left[
		\underline{Q}^{( i, j )}(\nu)
		,
		\underline{Q}^{( i, j )}(\nu)
		+
		Q^{( i, j )}(\nu)
	\right),
\end{align*}
where
\begin{align*}
	\underline{Q}^{( i, j )}(\nu)
	:&=
	\underset{ \mathfrak{i} \neq \mathfrak{j} }{
		\sum_{ \mathfrak{i} \leq i }
		\sum_{ \mathfrak{j} \leq \underline{j}( \mathfrak{i}; i ) }
	}
	Q^{( \mathfrak{i}, \mathfrak{j} )}(\nu),
	&
	&\mbox{with}	
	&
	\underline{j}( \mathfrak{i}; i )
	&=
	\begin{cases}
		n,		& \mathfrak{i} < i,
		\\
		j-1,	& \mathfrak{i} = i.
	\end{cases}
\end{align*}
Then, for a random environment $\mu$ and a Poisson random measure $N$ on $\RR_+ \times \RR_+$ with compensator $\diff t \otimes \diff r$, the process (or \textit{regime}) $\xi^\mu$, defined as
\begin{align}
    \label{Eq:RS-Y}
    \xi^\mu_t
    =
    \int_{ [0,t] \times \RR_+ } q(\mu_s,\xi^\mu_{s-},r) N(\diff s, \diff r),
\end{align}
with
\begin{align*}
    &q(\nu,i,r) := \sum_{j \in \Ee} (j-i) \Ind{ \Gamma^{(i,j)}(\nu) }(r),
    &
    &\forall (\nu,i,r) \in \Pp_2(\RR^d) \times \Ee \times \RR_+,
\end{align*}
is a Markov chain with conditional transition matrix $Q$ given $\mu$:
\begin{align*}
    \Prob{ \xi^\mu_{ t + \delta} = j \middle\vert \xi^\mu_t = i, \Ff^\mu_t  }
    =
    \begin{cases}
		Q^{(i,j)}(\mu_t)\delta + o(\delta),		& i \neq j,
		\\
		1 + Q^{(i,i)}(\mu_t)\delta + o(\delta),	& i=j.
    \end{cases}
\end{align*}

Observe that equation \eqref{Eq:RS-Y} allows us to express the regime $\eta^\mu$ in terms of an admissible common noise on $[0,T] \times \Ee$ by setting
\begin{align*}
    \xi_t^\mu
    =
    \int_{ \{t\} \times \Ee } e \ \eta^\lambda(\diff s, \diff e)
    =:
    \xi_{N^\lambda_t},
\end{align*}
where
\begin{align*}
    \eta^\lambda(\diff s, \diff e)
    =&
    \sum_{j \in \Ee} N\Parentheses{ \diff s, \Gamma^{(\xi_{s-}^\mu,j)}(\mu_{s-}) } 
        \otimes \delta_{j}(\diff e)
\end{align*}
is the lifting (see (\ref{lifting})) for the marked process $(N^\lambda,\xi)$, characterised by the admissible kernel
\begin{align*}
    K^\lambda(t,\diff e)
    =
    \sum_{j \in \Ee} Q^{(\xi^\mu_{t-},j)}(\mu_{t-})\ \delta_{j}(\diff e).
\end{align*}

\end{example}

There are many interesting and useful results that can be derived from the previous example, among them we  highlight the  role that  Markov chains of the type described above play in \textit{Regime-Switching McKean-Vlasov SDEs} and \textit{Regime-Switching Mean-Field Games}. In both cases, it is required to solve SDEs driven by an idiosyncratic noise process, usually a Brownian motion, where each of the coefficients involved in the dynamic is parametrized by a continuous time, irreducible Markov chain with finite state-space. An important difference with  the model presented here is that the Markov chain itself may depend on the law of the solution process. For example, from \cite{shao_propagation_2022} we know that the conditional McKean-Vlasov SDEs ($X$) with regime switching ($\xi$) satisfy the system of SDEs
\begin{align*}
    &\diff X_t = b(X_t, \Ll( X | \xi )_t , \xi_t ) \diff t + \sigma \diff W_t,
    \\
    &\diff \xi_t = \int_{\RR_+} q( \Ll( X | \xi )_t , \xi_t, r ) N(\diff s,\diff r).
\end{align*}
Fortunately, this fact  is captured by the notion of {\it admissibility} through the measure-valued process $\mu$, and so,  \textit{Common Noise} and \textit{Regime Switching} processes can be studied under the same framework.


\providecommand{\bysame}{\leavevmode\hbox to3em{\hrulefill}\thinspace}
\providecommand{\MR}{\relax\ifhmode\unskip\space\fi MR }
\providecommand{\MRhref}[2]{%
  \href{http://www.ams.org/mathscinet-getitem?mr=#1}{#2}
}



\begin{thebibliography}{10}

\bibitem{aksamit_enlargement_2017}
Anna Aksamit and Monique Jeanblanc, \emph{Enlargement of {{Filtration}} with
  {{Finance}} in {{View}}}, {{SpringerBriefs}} in {{Quantitative Finance}},
  {Springer International Publishing}, {Cham}, 2017.

\bibitem{baccelli_elements_2002}
Francois Baccelli and Pierre Bremaud, \emph{Elements of {{Queueing Theory}}:
  {{Palm Martingale Calculus}} and {{Stochastic Recurrences}}}, {Springer
  Science \& Business Media}, December 2002.

\bibitem{bacry_hawkes_2015}
Emmanuel Bacry, Iacopo Mastromatteo, and Jean-Fran{\c c}ois Muzy, \emph{Hawkes
  {{Processes}} in {{Finance}}}, Market Microstructure and Liquidity
  \textbf{01} (2015), no.~01, 1550005.

\bibitem{billingsley_convergence_2013}
Patrick Billingsley, \emph{Convergence of {{Probability Measures}}}, {John
  Wiley \& Sons}, June 2013.

\bibitem{bismut_conjugate_1973}
Jean-Michel Bismut, \emph{Conjugate convex functions in optimal stochastic
  control}, Journal of Mathematical Analysis and Applications \textbf{44}
  (1973), no.~2, 384--404.

\bibitem{bismut_introductory_1978}
\bysame, \emph{An {{Introductory Approach}} to {{Duality}} in {{Optimal
  Stochastic Control}}}, SIAM Review \textbf{20} (1978), no.~1, 62--78.

\bibitem{bogachev_measure_2007}
Vladimir~I. Bogachev, \emph{Measure {{Theory}}}, {Springer Berlin Heidelberg},
  {Berlin, Heidelberg}, 2007.

\bibitem{bremaud_point_2020}
Pierre Br{\`e}maud, \emph{Point {{Process Calculus}} in {{Time}} and {{Space}}:
  {{An Introduction}} with {{Applications}}}, Probability {{Theory}} and
  {{Stochastic Modelling}}, vol.~98, {Springer International Publishing},
  {Cham}, 2020.

\bibitem{bremaud_changes_1978}
Pierre Br{\`e}maud and Marc Yor, \emph{Changes of filtrations and of
  probability measures}, Zeitschrift f\"ur Wahrscheinlichkeitstheorie und
  Verwandte Gebiete \textbf{45} (1978), no.~4, 269--295.

\bibitem{caines_large_2006}
Peter~E. Caines, Minyi Huang, and Roland~P. Malham{\'e}, \emph{Large population
  stochastic dynamic games: Closed-loop {{McKean-Vlasov}} systems and the
  {{Nash}} certainty equivalence principle}, Communications in Information and
  Systems \textbf{6} (2006), no.~3, 221--252.

\bibitem{cardaliaguet_notes_nodate}
Pierre Cardaliaguet, \emph{Notes on {{Mean Field Games}}}, Tech. report.

\bibitem{cardaliaguet_master_2019}
\bysame, \emph{The master equation and the convergence problem in mean field
  games}, Annals of Mathematics Studies, no. 201, {Princeton University Press},
  {Princeton, NJ}, 2019.

\bibitem{carmona_probabilistic_2013}
Ren{\'e} Carmona and Fran{\c c}ois Delarue, \emph{Probabilistic {{Analysis}} of
  {{Mean-Field Games}}}, SIAM Journal on Control and Optimization \textbf{51}
  (2013), no.~4, 2705--2734.

\bibitem{carmona_probabilistic_2018}
\bysame, \emph{Probabilistic {{Theory}} of {{Mean Field Games}} with
  {{Applications I}}}, Probability {{Theory}} and {{Stochastic Modelling}},
  vol.~83, {Springer International Publishing}, {Cham}, 2018.

\bibitem{carmona_probabilistic_2018-1}
\bysame, \emph{Probabilistic {{Theory}} of {{Mean Field Games}} with
  {{Applications II}}}, Probability {{Theory}} and {{Stochastic Modelling}},
  vol.~84, {Springer International Publishing}, {Cham}, 2018.

\bibitem{carmona_mean_2016}
Ren{\'e} Carmona, Fran{\c c}ois Delarue, and Daniel Lacker, \emph{Mean field
  games with common noise}, The Annals of Probability \textbf{44} (2016),
  no.~6.

\bibitem{chaintron_propagation_2022}
Louis-Pierre Chaintron and Antoine Diez, \emph{Propagation of chaos: A review
  of models, methods and applications. {{I}}. {{Models}} and methods}, May
  2022.

\bibitem{chaintron_propagation_2022-1}
\bysame, \emph{Propagation of chaos: A review of models, methods and
  applications. {{II}}. {{Applications}}}, May 2022.

\bibitem{cohen_stochastic_2015}
Samuel~N. Cohen and Robert~James Elliott, \emph{Stochastic calculus and
  applications}, 2nd ed ed., Probability and {{Its Applications}},
  {Birkh\"auser}, {New York, NY}, 2015.

\bibitem{daley_introduction_2003}
D.~J. Daley and D.~{Vere-Jones}, \emph{An {{Introduction}} to the {{Theory}} of
  {{Point Processes}}. {{Volume I}}: {{Elementary Theory}} and {{Methods}}},
  Probability and Its {{Applications}}, {Springer-Verlag}, {New York}, 2003.

\bibitem{daley_introduction_2008}
\bysame, \emph{An {{Introduction}} to the {{Theory}} of {{Point Processes}}.
  {{Volume II}}: {{General Theory}} and {{Structure}}}, Probability and {{Its
  Applications}}, {Springer}, {New York, NY}, 2008.

\bibitem{el_karoui_backward_1997}
N.~El~Karoui, S.~Peng, and M.~C. Quenez, \emph{Backward {{Stochastic
  Differential Equations}} in {{Finance}}}, Mathematical Finance \textbf{7}
  (1997), no.~1, 1--71.

\bibitem{azema_vershiks_2001}
Michel {\'E}mery and Walter Schachermayer, \emph{On {{Vershik}}'s
  {{Standardness Criterion}} and {{Tsirelson}}'s {{Notion}} of {{Cosiness}}},
  S\'eminaire de {{Probabilit\'es XXXV}} (J.~Az{\'e}ma, M.~{\'E}mery,
  M.~Ledoux, and M.~Yor, eds.), vol. 1755, {Springer Berlin Heidelberg},
  {Berlin, Heidelberg}, 2001, pp.~265--305.

\bibitem{fouque_handbook_2013}
Jean-Pierre Fouque and Joseph~A. Langsam, \emph{Handbook on {{Systemic Risk}}},
  {Cambridge University Press}, May 2013.

\bibitem{hawkes_hawkes_2018}
Alan~G. Hawkes, \emph{Hawkes processes and their applications to finance: A
  review}, Quantitative Finance \textbf{18} (2018), no.~2, 193--198.

\bibitem{hawkes_hawkes_2022}
\bysame, \emph{Hawkes jump-diffusions and finance: A brief history and review},
  The European Journal of Finance \textbf{28} (2022), no.~7, 627--641.

\bibitem{jacod_calcul_2006}
J.~Jacod, \emph{{Calcul Stochastique et Probl\`emes de Martingales}},
  {Springer}, November 2006.

\bibitem{jacod_limit_2003}
Jean Jacod and Albert~N. Shiryaev, \emph{Limit {{Theorems}} for {{Stochastic
  Processes}}}, Grundlehren Der Mathematischen {{Wissenschaften}}, vol. 288,
  {Springer Berlin Heidelberg}, {Berlin, Heidelberg}, 2003.

\bibitem{kurtz_weak_2014}
Thomas Kurtz, \emph{Weak and strong solutions of general stochastic models},
  Electronic Communications in Probability \textbf{19} (2014), no.~none.

\bibitem{ma_forward-backward_2007}
Jin Ma and Jiongmin Yong, \emph{Forward-{{Backward Stochastic Differential
  Equations}} and their {{Applications}}}, Lecture {{Notes}} in
  {{Mathematics}}, vol. 1702, {Springer Berlin Heidelberg}, {Berlin,
  Heidelberg}, 2007.

\bibitem{oettli_existence_1998}
W.~Oettli and D.~Schl{\"a}ger, \emph{Existence of {{Equilibria}} for
  {{G-Monotone Mappings}}}, IFAC Proceedings Volumes \textbf{31} (1998),
  no.~13, 13--17.

\bibitem{brokate_generalized_nodate}
Werner Oettli and Dirk Schl{\"a}ger, \emph{Generalized vectorial equilibria and
  generalized monotonicity}, Functional {{Analysis}} with {{Current
  Applications}} in {{Science}}, {{Technology}} and {{Industry}} (Martin
  Brokate and Abul~Hasan Siddiqi, eds.), {Chapman and Hall/CRC}, zeroth ed.

\bibitem{oksendal_applied_2019}
Bernt {\O}ksendal and Agn{\`e}s Sulem, \emph{Applied {{Stochastic Control}} of
  {{Jump Diffusions}}}, Universitext, {Springer International Publishing},
  {Cham}, 2019.

\bibitem{papapantoleon_existence_2018}
Antonis Papapantoleon, Dylan Possama{\"i}, and Alexandros Saplaouras,
  \emph{Existence and uniqueness results for {{BSDE}} with jumps: The whole
  nine yards}, Electronic Journal of Probability \textbf{23} (2018), no.~none,
  1--68.

\bibitem{pardoux_adapted_1990}
E.~Pardoux and S.~G. Peng, \emph{Adapted solution of a backward stochastic
  differential equation}, Systems \& Control Letters \textbf{14} (1990), no.~1,
  55--61.

\bibitem{peng_general_1990}
Shige Peng, \emph{A {{General Stochastic Maximum Principle}} for {{Optimal
  Control Problems}}}, SIAM Journal on Control and Optimization \textbf{28}
  (1990), no.~4, 966--979.

\bibitem{peng_fully_1999}
Shige Peng and Zhen Wu, \emph{Fully {{Coupled Forward-Backward Stochastic
  Differential Equations}} and {{Applications}} to {{Optimal Control}}}, SIAM
  Journal on Control and Optimization \textbf{37} (1999), no.~3, 825--843.

\bibitem{protter_stochastic_2005}
Philip~E. Protter, \emph{Stochastic {{Integration}} and {{Differential
  Equations}}}, Stochastic {{Modelling}} and {{Applied Probability}}, vol.~21,
  {Springer Berlin Heidelberg}, {Berlin, Heidelberg}, 2005.

\bibitem{shao_propagation_2022}
Jinghai Shao and Dong Wei, \emph{Propagation of chaos and conditional
  {{McKean-Vlasov SDEs}} with regime-switching}, Frontiers of Mathematics in
  China \textbf{17} (2022), no.~4, 731--746.

\bibitem{burkholder_topics_1991}
Alain-Sol Sznitman, \emph{Topics in propagation of chaos}, vol. 1464,
  pp.~165--251, {Springer Berlin Heidelberg}, {Berlin, Heidelberg}, 1991.

\bibitem{tang_necessary_1994}
Shanjian Tang and Xunjing Li, \emph{Necessary {{Conditions}} for {{Optimal
  Control}} of {{Stochastic Systems}} with {{Random Jumps}}}, SIAM Journal on
  Control and Optimization \textbf{32} (1994), no.~5, 1447--1475.

\bibitem{yin_hybrid_2010}
G.~George Yin and Chao Zhu, \emph{Hybrid {{Switching Diffusions}}}, Stochastic
  {{Modelling}} and {{Applied Probability}}, vol.~63, {Springer New York}, {New
  York, NY}, 2010.

\bibitem{yong_stochastic_1999}
J.~Yong and Xun~Yu Zhou, \emph{Stochastic controls: {{Hamiltonian}} systems and
  {{HJB}} equations}, Applications of Mathematics, no.~43, {Springer}, {New
  York}, 1999.

\bibitem{yong_finding_1997}
Jiongmin Yong, \emph{Finding adapted solutions of forward\textendash backward
  stochastic differential equations: Method of continuation}, Probability
  Theory and Related Fields \textbf{107} (1997), no.~4, 537--572.

\bibitem{zhen_forward-backward_1999}
Wu~Zhen, \emph{Forward-backward stochastic differential equations with
  {{Brownian}} motion and poisson process}, Acta Mathematicae Applicatae Sinica
  \textbf{15} (1999), no.~4, 433--443.

\end{thebibliography}
\end{document}